\documentclass[11pt]{article}%
\usepackage{amsmath}
\usepackage{amsfonts}
\usepackage{amssymb}
\usepackage{graphicx}
\usepackage{theorem}
\usepackage{makeidx}
\usepackage{verbatim}
\usepackage[colorlinks=true]{hyperref}%
\providecommand{\U}[1]{\protect\rule{.1in}{.1in}}
\newtheorem{thm}{Theorem}[section]
\newtheorem{lm}[thm]{Lemma}
\newtheorem{pr}[thm]{Proposition}
\newtheorem{df}[thm]{Definition}

\newtheorem{cor}[thm]{Corollary}
{\theorembodyfont{\upshape}
\newtheorem{examp}[thm]{Example}
}
\numberwithin{equation}{section} 
\begin{document}

\title{ }

\begin{center}
\vspace*{1.5cm}

\textbf{OPENNESS STABILITY AND IMPLICIT MULTIFUNCTION THEOREMS.}

\textbf{APPLICATIONS\ TO\ VARIATIONAL SYSTEMS}

\vspace*{1cm}

M. DUREA

{\small {Faculty of Mathematics, \textquotedblleft Al. I. Cuza" University,} }

{\small {Bd. Carol I, nr. 11, 700506 -- Ia\c{s}i, Romania,} }

{\small {e-mail: \texttt{durea@uaic.ro}}}

\bigskip

R. STRUGARIU

{\small Department of Mathematics, \textquotedblleft Gh. Asachi" Technical
University, }

{\small {Bd. Carol I, nr. 11, 700506 -- Ia\c{s}i, Romania,} }

{\small {e-mail: \texttt{rstrugariu@tuiasi.ro}}}
\end{center}

\bigskip

\bigskip

\noindent{\small {\textbf{Abstract: }}In this paper we aim to present two
general results regarding, on one hand, the openness stability of set-valued
maps and, on the other hand, the metric regularity behavior of the implicit
multifunction related to a generalized variational system. Then, these results
are applied in order to obtain, in a natural way, and in a widely studied
case, several relations between the metric regularity moduli of the field maps
defining the variational system and the solution map. Our approach allows us
to complete and extend several very recent results in literature.}

\bigskip


\bigskip


\section{Introduction}

This paper belongs to the active area of research concerning parametric
variational systems and it aims to enter into dialog with some very recent
works of Arag\'{o}n Artacho and Mordukhovich (\cite{ArtMord2009},
\cite{ArtMord2010}) and Ngai, Tron and Th\'{e}ra (\cite{NTT}). Note that, in
turn, these papers extend many results of Dontchev and Rockafellar
(\cite{DontRock2009a}, \cite{DontRock2009b}).

Firstly, our research on the relations between metric regularity/Lipschitz
moduli of an initial parametric field map and the associated implicit
multifunction map led us to the rediscovery of a very nice Theorem of Ursescu
\cite{Urs1996} concerning the stability of openness of set-valued maps
(Theorem \ref{main} below). The proof we provide here for this result is
appropriate enough for getting some extra assertions compared with the initial
paper of Ursescu. Later in the paper, this result is a key ingredient in order
to get a natural and precise answer to the question of how regularity
constants of the involved maps relate each other.

Secondly, we were interested in enlarging the framework commonly used as being
the defining form of a parametric variational system to the case of a general
field map. To be more specific, let $X,Y,P$ be Banach spaces, $H:X\times
P\rightrightarrows Y$ be a multifunction and define the implicit set-valued
map $S:P\rightrightarrows X$ by%
\[
S(p)=\{x\in X\mid0\in H(x,p)\}.
\]
Then, in this second part of our work we find how metric regularity and
Lipschitz properties of $H$ and $S$ are related under certain assumptions
(Theorem \ref{impl}). This second main result is in fact a general implicit
multifunction theorem. After fixing these two main tools we are able to
present in a natural way the situation where $H$ is a sum of two set-valued
maps $F,G$ of the form $H(x,p)=F(x,p)+G(x)$. Note that this case is more
general than the situation considered in \cite{ArtMord2009} by the presence of
the set-valued map $F$ instead of a single-valued map. Moreover, this is the
most general situation one can consider because it is not possible to get good
results concerning the Lipschitz properties when $G$ depends on the parameter
$p,$ in virtue of \cite[Remark 3.6. (iii)]{ArtMord2009}. However, note that in
\cite{NTT} the authors deal with a sort of metric regularity of the solution
map associated to the sum of two parametric set-valued maps. In our framework,
when we put at work together the two main results, we are in the position to
indicate in a smooth manner the relations between the regularity moduli of
$S,G$ and $F.$ We hope that our main results and their combination will bring
more light on the previous results on this topic.

We would like to mention that, in comparison with \cite{ArtMord2009} and
\cite{ArtMord2010}, we get here only results in which the assumptions are on
$F$ and $G$ and the conclusion concerns $S.$ The converse situation considered
in the quoted works (from $S$ to $G$) is not presented here because it has too
many similarities with the corresponding results of Artacho and Mordukhovich.
Any interested reader could find the arguments to obtain such results in our
framework, but with Artacho and Mordukhovich tools.

The paper is organized as follows. In the next section we present the
notations, the concepts and the basic facts we use in the sequel. The third
section contains the main results of the paper we have presented in few words
above. The last section investigates the widely studied form of the parametric
variational systems we can find in literature. We show here how the main
results concerning the stability of the linear openness of set-valued maps and
the implicit multifunctions could be combined in order to get quite easily the
estimation of the regularity constants of the solution map.

\section{Preliminaries}

This section contains some basic definitions and results used in the sequel.
In what follows, we suppose that all the involved spaces are Banach. In this
setting, $B(x,r)$ and $D(x,r)$ denote the open and the closed ball with center
$x$ and radius $r,$ respectively$.$ Sometimes we write $\mathbb{D}_{X}$ for
the closed unit ball of $X$. If $x\in X$ and $A\subset X,$ one defines the
distance from $x$ to $A$ as $d(x,A):=\inf\{\left\Vert x-a\right\Vert \mid a\in
A\}.$ As usual, we use the convention $d(x,\emptyset)=\infty.$ For a non-empty
set $A\subset X$ we put $\operatorname*{cl}A,$ $\operatorname*{int}A$ for the
topological closure and interior, respectively. Also, a set $A$ is said to be
locally closed if for every $a\in A,$ there exists $r>0$ such that the set
$A\cap D(a,r)$ is closed. If $\overline{a}\in A$ and $A\cap V$ is locally
closed, where $V$ is a neighborhood of $\overline{a},$ we say that the set $A$
is locally closed around $\overline{a}.$

Consider now a multifunction $F:X\rightrightarrows Y$. The domain and the
graph of $F$ are denoted respectively by
\[
\operatorname*{Dom}F:=\{x\in X\mid F(x)\neq\emptyset\}
\]
and%
\[
\operatorname*{Gr}F=\{(x,y)\in X\times Y\mid y\in F(x)\}.
\]
If $A\subset X$ then $F(A):=%
{\displaystyle\bigcup\limits_{x\in A}}
F(x).$ The inverse set-valued map of $F$ is $F^{-1}:Y\rightrightarrows X$
given by $F^{-1}(y)=\{x\in X\mid y\in F(x)\}$.

Recall that a multifunction $F$ is inner semicontinuous at $(x,y)\in
\operatorname*{Gr}F$ if for every open set $D\subset Y$ with $y\in D,$ there
exists a neighborhood $U\in\mathcal{V}(x)$ such that for every $x^{\prime}\in
U,$ $F(x^{\prime})\cap D\neq\emptyset$ (where $\mathcal{V}(x)$ stands for the
system of the neighborhoods of $x$).

We remind now the concepts of openness at linear rate, metric regularity and
Lipschitz-likeness of a multifunction around the reference point.

\begin{df}
\label{around}Let $L>0,$ $F:X\rightrightarrows Y$ be a multifunction and
$(\overline{x},\overline{y})\in\operatorname{Gr}F.$

(i) $F$ is said to be open at linear rate $L>0,$ or $L-$open around
$(\overline{x},\overline{y})$ if there exist a positive number $\varepsilon>0$
and two neighborhoods $U\in\mathcal{V}(\overline{x}),$ $V\in\mathcal{V}%
(\overline{y})$ such that, for every $\rho\in]0,\varepsilon\lbrack$ and every
$(x,y)\in\operatorname*{Gr}F\cap\lbrack U\times V],$%
\begin{equation}
B(y,\rho L)\subset F(B(x,\rho)). \label{Lopen}%
\end{equation}

The supremum of $L>0$ over all the combinations $(L,U,V,\varepsilon)$ for
which (\ref{Lopen}) holds is denoted by $\operatorname*{lop}F(\overline
{x},\overline{y})$ and is called the exact linear openness bound, or the exact
covering bound of $F$ around $(\overline{x},\overline{y}).$

(ii) $F$ is said to be Lipschitz-like, or has Aubin property around
$(\overline{x},\overline{y})$ with constant $L>0$ if there exist two
neighborhoods $U\in\mathcal{V}(\overline{x}),$ $V\in\mathcal{V}(\overline{y})$
such that, for every $x,u\in U,$%
\begin{equation}
F(x)\cap V\subset F(u)+L\left\Vert x-u\right\Vert \mathbb{D}_{Y}.
\label{LLip_like}%
\end{equation}

The infimum of $L>0$ over all the combinations $(L,U,V)$ for which
(\ref{LLip_like}) holds is denoted by $\operatorname*{lip}F(\overline
{x},\overline{y})$ and is called the exact Lipschitz bound of $F$ around
$(\overline{x},\overline{y}).$

(iii) $F$ is said to be metrically regular around $(\overline{x},\overline
{y})$ with constant $L>0$ if there exist two neighborhoods $U\in
\mathcal{V}(\overline{x}),$ $V\in\mathcal{V}(\overline{y})$ such that, for
every $(x,y)\in U\times V,$%
\begin{equation}
d(x,F^{-1}(y))\leq Ld(y,F(x)). \label{Lmreg}%
\end{equation}

The infimum of $L>0$ over all the combinations $(L,U,V)$ for which
(\ref{Lmreg}) holds is denoted by $\operatorname*{reg}F(\overline{x}%
,\overline{y})$ and is called the exact regularity bound of $F$ around
$(\overline{x},\overline{y}).$
\end{df}

The next proposition contains the well-known links between the notions
presented above. See \cite[Theorems 1.49, 1.52]{Mor2006} for more details
about the proof.

\begin{pr}
\label{link_around}Let $F:X\rightrightarrows Y$ be a multifunction and
$(\overline{x},\overline{y})\in\operatorname{Gr}F.$ Then $F$ is open at linear
rate around $(\overline{x},\overline{y})$ iff $F^{-1}$ is Lipschitz-like
around $(\overline{y},\overline{x})$ iff $F$ is metrically regular around
$(\overline{x},\overline{y})$. Moreover, in every of the previous situations,%
\[
(\operatorname*{lop}F(\overline{x},\overline{y}))^{-1}=\operatorname*{lip}%
F^{-1}(\overline{y},\overline{x})=\operatorname*{reg}F(\overline{x}%
,\overline{y}).
\]

\end{pr}

It is well known that the corresponding \textquotedblleft at point" properties
are significantly different from the \textquotedblleft around point" ones. Let
us introduce now some of these notions. For more related concepts we refer to
\cite{ArtMord2009}.

\begin{df}
\label{at}Let $L>0,$ $F:X\rightrightarrows Y$ be a multifunction and
$(\overline{x},\overline{y})\in\operatorname{Gr}F.$

(i) $F$ is said to be open at linear rate $L,$ or $L-$open at $(\overline
{x},\overline{y})$ if there exists a positive number $\varepsilon>0$ such
that, for every $\rho\in]0,\varepsilon\lbrack,$%
\begin{equation}
B(\overline{y},\rho L)\subset F(B(\overline{x},\rho)). \label{Lopen_at}%
\end{equation}

The supremum of $L>0$ over all the combinations $(L,\varepsilon)$ for which
(\ref{Lopen_at}) holds is denoted by $\operatorname*{plop}F(\overline
{x},\overline{y})$ and is called the exact punctual linear openness bound of
$F$ at $(\overline{x},\overline{y}).$

(ii) $F$ is said to be pseudocalm with constant $L,$ or $L-$pseudocalm at
$(\overline{x},\overline{y}),$ if there exists a neighborhood $U\in
\mathcal{V}(\overline{x})$ such that, for every $x\in U,$%
\begin{equation}
d(\overline{y},F(x))\leq L\left\Vert x-\overline{x}\right\Vert .
\label{Lpseudocalm}%
\end{equation}

The infimum of $L>0$ over all the combinations $(L,U)$ for which
(\ref{Lpseudocalm}) holds is denoted by $\operatorname*{psdclm}F(\overline
{x},\overline{y})$ and is called the exact bound of pseudocalmness for $F$ at
$(\overline{x},\overline{y}).$

(iii) $F$ is said to be metrically hemiregular with constant $L,$ or
$L-$metrically hemiregular at $(\overline{x},\overline{y})$ if there exists a
neighborhood $V\in\mathcal{V}({\overline{y}})$ such that, for every $y\in V,$%
\begin{equation}
d(\overline{x},F^{-1}(y))\leq L\left\Vert y-{\overline{y}}\right\Vert .
\label{Lhemireg}%
\end{equation}

The infimum of $L>0$ over all the combinations $(L,V)$ for which
(\ref{Lhemireg}) holds is denoted by $\operatorname*{hemreg}F(\overline
{x},\overline{y})$ and is called the exact hemiregularity bound of $F$ at
$(\overline{x},\overline{y}).$
\end{df}

The term of metric hemiregularity appears in \cite[Definition 5.1]%
{ArtMord2010}, where the link with \textquotedblleft Lipschitz lower
semicontinuity" (i.e., pseudocalmness in our terminology) of the inverse
multifunction is emphasized. The notion of pseudocalmness is used under the
term of $L-$Lipschitz in \cite{Urs1996}, where other concepts of relative
openness and relative $L-$Lipschitz properties are introduced and discussed.

The next proposition lists some equivalences between these \textquotedblleft
at point" notions. We give the (elementary) proof for the completeness.

\begin{pr}
\label{link_at}Let $L>0,$ $F:X\rightrightarrows Y$ and $(\overline
{x},\overline{y})\in\operatorname{Gr}F.$ Then $F$ is $L-$open at
$(\overline{x},\overline{y})$ iff $F^{-1}$ is $L^{-1}-$pseudocalm at
$(\overline{y},\overline{x})$ iff $F$ is $L^{-1}-$metrically hemiregular at
$(\overline{x},\overline{y})$. Moreover, in every of the previous situations,%
\[
(\operatorname*{plop}F(\overline{x},{\overline{y}}))^{-1}%
=\operatorname*{psdclm}F^{-1}({\overline{y},}\overline{x}%
)=\operatorname*{hemreg}F(\overline{x},{\overline{y}}).
\]

\end{pr}

\noindent\textbf{Proof.} It's obvious from the very definitions that $F^{-1}$
is $L^{-1}-$pseudocalm at $(\overline{y},\overline{x})$ iff $F$ is $L^{-1}%
-$metrically hemiregular at $(\overline{x},\overline{y})$. Suppose now that
$F$ is $L-$open at $(\overline{x},\overline{y}).$ Then there exists
$\varepsilon>0$ such that, for every $\rho\in]0,\varepsilon\lbrack,$
(\ref{Lopen_at}) holds. Consider $\varepsilon^{\prime}:=L\varepsilon$ and take
arbitrarily $y\in B(\overline{y},\varepsilon^{\prime}).$ Then there exist
$\rho\in]0,\varepsilon\lbrack$ and $\gamma$ arbitrary small such that
$\left\Vert y-\overline{y}\right\Vert =L\rho<L(\rho+\gamma)<L\varepsilon.$
Using the $L-$openness of $F$ at $(\overline{x},\overline{y}),$ $y\in
B(\overline{y},L(\rho+\gamma))\subset F(B(\overline{x},\rho+\gamma)).$
Consequently, one can find $x\in B(\overline{x},\rho+\gamma)\cap F^{-1}(y)$.
Then, for $\gamma$ arbitrary small, $d(\overline{x},F^{-1}(y))\leq\left\Vert
x-\overline{x}\right\Vert <\rho+\gamma,$ whence $d(\overline{x},F^{-1}%
(y))\leq\rho=L^{-1}\left\Vert y-\overline{y}\right\Vert ,$ and the first
implication is now proved.

Suppose now that $F^{-1}$ is $L^{-1}-$pseudocalm at $(\overline{y}%
,\overline{x}),$ so there exists $\varepsilon>0$ such that, for every $y\in
B(\overline{y},\varepsilon),$ $d(\overline{x},F^{-1}(y))\leq L^{-1}\left\Vert
y-\overline{y}\right\Vert .$ Take $\varepsilon^{\prime}:=L^{-1}\varepsilon$
and arbitrarily $\rho\in]0,\varepsilon^{\prime}[.$ If $y^{\prime}\in
B(\overline{y},L\rho),$ we obtain that $\left\Vert y^{\prime}-\overline
{y}\right\Vert <L\varepsilon^{\prime}=\varepsilon,$ so $d(\overline{x}%
,F^{-1}(y^{\prime}))\leq L^{-1}\left\Vert y^{\prime}-\overline{y}\right\Vert
<\rho.$ Consequently, there exists $x^{\prime}\in F^{-1}(y^{\prime})$ such
that $\left\Vert x^{\prime}-\overline{x}\right\Vert <\rho,$ whence $y^{\prime
}\in F(x^{\prime})\subset F(B({\overline{x}},\rho)).$ The proof is now
complete.$\hfill\square$

\bigskip

See \cite[Section 11]{Urs2008} for an example of a multifunction which is open
at linear rate at a point, hence on the basis of Proposition \ref{link_at} is
metrically hemiregular at this point, even if that multifunction is not
metrically regular around any point.

Recall that $\mathcal{L}(X,Y)$ denotes the normed vector space of linear
bounded operators acting between $X$ and $Y$. If $A\in\mathcal{L}(X,Y),$ then
the \textquotedblleft at" and \textquotedblleft around point" notions do
coincide. In fact, $A$ is metrically regular around every $x\in X$ iff $A$ is
metrically hemiregular at every $x\in X$ iff $A$ is open with linear rate
around every $x\in X$ iff $A$ is open with linear rate at every $x\in X$ iff
$A$ is surjective. Moreover, in every of these cases we have%
\[
\operatorname*{hemreg}A=\operatorname*{reg}A=(\operatorname*{plop}%
A)^{-1}=(\operatorname*{lop}A)^{-1}=\left\Vert (A^{\ast})^{-1}\right\Vert ,
\]

\noindent where $A^{\ast}\in\mathcal{L}(Y^{\ast},X^{\ast})$ denotes the
adjoint operator and $\operatorname*{hemreg}A,$ $\operatorname*{reg}A,$
$\operatorname*{plop}A$ and $\operatorname*{lop}A$ are common for all the
points $x\in X$ (see, for more details, \cite[Proposition 5.2]{ArtMord2010}).

\bigskip

In the following, we introduce the corresponding partial notions of linear
openness, metric regularity and Lipschitz-like property around the reference
point for a parametric set-valued map. Below, we denote by $P$ the Banach
space of parameters.

\begin{df}
Let $L>0,$ $F:X\times P\rightrightarrows Y$ be a multifunction, $((\overline
{x},\overline{p}),\overline{y})\in\operatorname{Gr}F$ and for every $p\in P,$
denote $F_{p}(\cdot):=F(\cdot,p).$

(i) $F$ is said to be open at linear rate $L>0,$ or $L-$open, with respect to
$x$ uniformly in $p$ around $((\overline{x},\overline{p}),\overline{y})$ if
there exist a positive number $\varepsilon>0$ and some neighborhoods
$U\in\mathcal{V}(\overline{x}),$ $V\in\mathcal{V}(\overline{p}),$
$W\in\mathcal{V}(\overline{y})$ such that, for every $\rho\in]0,\varepsilon
\lbrack,$ every $p\in V$ and every $(x,y)\in\operatorname*{Gr}F_{p}\cap\lbrack
U\times W],$%
\begin{equation}
B(y,\rho L)\subset F_{p}(B(x,\rho)), \label{pLopen}%
\end{equation}

The supremum of $L>0$ over all the combinations $(L,U,V,W,\varepsilon)$ for
which (\ref{pLopen}) holds is denoted by $\widehat{\operatorname*{lop}}%
_{x}F((\overline{x},\overline{p}),\overline{y})$ and is called the exact
linear openness bound, or the exact covering bound of $F$ in $x$ around
$((\overline{x},\overline{p}),\overline{y}).$

(ii) $F$ is said to be Lipschitz-like, or has Aubin property, with respect to
$x$ uniformly in $p$ around $((\overline{x},\overline{p}),\overline{y})$ with
constant $L>0$ if there exist some neighborhoods $U\in\mathcal{V}(\overline
{x}),$ $V\in\mathcal{V}(\overline{p}),$ $W\in\mathcal{V}(\overline{y})$ such
that, for every $x,u\in U$ and every $p\in V,$%
\begin{equation}
F_{p}(x)\cap W\subset F_{p}(u)+L\left\Vert x-u\right\Vert \mathbb{D}_{Y}.
\label{pLLip_like}%
\end{equation}

The infimum of $L>0$ over all the combinations $(L,U,V,W)$ for which
(\ref{pLLip_like}) holds is denoted by $\widehat{\operatorname*{lip}}%
_{x}F((\overline{x},\overline{p}),\overline{y})$ and is called the exact
Lipschitz bound of $F$ in $x$ around $((\overline{x},\overline{p}%
),\overline{y}).$

(iii) $F$ is said to be metrically regular with respect to $x$ uniformly in
$p$ around $((\overline{x},\overline{p}),\overline{y})$ with constant $L>0$ if
there exist some neighborhoods $U\in\mathcal{V}(\overline{x}),$ $V\in
\mathcal{V}(\overline{p}),$ $W\in\mathcal{V}(\overline{y})$ such that, for
every $(x,p,y)\in U\times V\times W,$%
\begin{equation}
d(x,F_{p}^{-1}(y))\leq Ld(y,F_{p}(x)). \label{pLmreg}%
\end{equation}

The infimum of $L>0$ over all the combinations $(L,U,V,W)$ for which
(\ref{pLmreg}) holds is denoted by $\widehat{\operatorname*{reg}}%
_{x}F((\overline{x},\overline{p}),\overline{y})$ and is called the exact
regularity bound of $F$ in $x$ around $((\overline{x},\overline{p}%
),\overline{y}).$
\end{df}

Similarly, one can define the notions of linear openness, metric regularity
and Lipschitz-like property with respect to $p$ uniformly in $x,$ and the
corresponding exact bounds.

\section{Main results}

We begin our analysis with an interesting result due to Ursescu (see, e.g.,
\cite[Theorem 1]{Urs1996}), which brings into the light the key fact that the
linear openness property of a difference type multifunction can be deduced
from the corresponding linear openness properties of its terms. This result
can be viewed as a deep generalization of the Graves Theorem (see the remark
after Theorem \ref{main}). Different variants of this result were stated in
\cite{Ioffe2000}, \cite{Dmi2005}, and the more general assumptions are that
$Y$ is a linear metric space with shift invariant metric, $X$ is metric space,
and the graphs of the involved multifunctions are locally complete. Moreover,
its importance seems to be crucial, because it can be putted into relation
with a large number of (very) actual topics, as the strongly regular
generalized equations of Robinson type (see \cite{Rob1980}), allowing for the
first time to deal with problems where both the terms are multivalued, but
also with the inverse and implicit type theorems for set-valued mappings, with
parametric variational inclusions and more, as one can see next. The proof we
provide here is somehow similar with that of \cite[Theorem 5]{Urs1996}, but it
is more direct and simple. Besides the result itself, this new proof of it is
a key element in extending the framework to the parametric case and,
furthermore, to the rest of the applications we make precise in the last section.

\begin{thm}
\label{main}Let $F:X\rightrightarrows Y$ and $G:Y\rightrightarrows X$ be two
multifunctions such that $\operatorname*{Gr}F$ and $\operatorname*{Gr}G$ are
locally closed. Suppose that $\operatorname*{Dom}(F-G^{-1})$ and
$\operatorname*{Dom}(G-F^{-1})$\ are nonempty and let $L>0$ and $M>0$ be such
that $LM>1.$ If $F$ is $L-$open at every point of its graph, and $G$ is
$M-$open at every point of its graph, then $F-G^{-1}$ is $(L-M^{-1})-$open at
every point of its graph and $G-F^{-1}$ is $(M-L^{-1})-$open at every point of
its graph.
\end{thm}

\noindent\textbf{Proof.} We prove only the first assertion, the other one
being completely symmetrical. Let $(x,w)\in\operatorname*{Gr}(F-G^{-1}).$ Then
there exist $y\in F(x)$ and $z\in G^{-1}(x)$ such that $w=y-z.$ Define the
multifunction $(F,G^{-1}):X\rightrightarrows Y\times Y$ by $(F,G^{-1}%
)(x):=F(x)\times G^{-1}(x)$ and remark that $(x,y,z)\in\operatorname*{Gr}%
(F,G^{-1}).$ Because $\operatorname*{Gr}F$ and $\operatorname*{Gr}G$ are
locally closed, it follows that $\operatorname*{Gr}(F,G^{-1})$ is locally
closed and one can find $\rho>0$ such that $\operatorname*{Gr}(F,G^{-1}%
)\cap\operatorname*{cl}W$ is closed, where%
\begin{equation}
W:=B(x,\rho)\times B(y,L\rho)\times B(z,M^{-1}\rho). \label{W}%
\end{equation}

Take $u\in B(w,(L-M^{-1})\rho).$ We must prove that $u\in(F-G^{-1}%
)(B(x,\rho)).$ One can find $\tau\in]0,1[$ such that $\left\Vert
u-w\right\Vert <\tau(L-M^{-1})\rho.$ Endow the space $X\times Y\times Y$ with
the norm%
\[
\left\Vert (p,q,r)\right\Vert _{0}:=\tau(L-M^{-1})\max\{\left\Vert
p\right\Vert ,L^{-1}\left\Vert q\right\Vert ,M\left\Vert r\right\Vert \}
\]

\noindent and apply the Ekeland variational principle to the function
$h:\operatorname{Gr}(F,G^{-1})\cap\operatorname*{cl}W\rightarrow\mathbb{R}%
_{+},$%
\[
h(p,q,r):=\left\Vert u-(q-r)\right\Vert .
\]

Then one can find a point $(a,b,c)\in\operatorname{Gr}(F,G^{-1})\cap
\operatorname*{cl}W$ such that%
\begin{equation}
\left\Vert u-(b-c)\right\Vert \leq\left\Vert u-(y-z)\right\Vert -\left\Vert
(a,b,c)-(x,y,z)\right\Vert _{0} \label{ek1}%
\end{equation}

\noindent and%
\begin{equation}
\left\Vert u-(b-c)\right\Vert \leq\left\Vert u-(q-r)\right\Vert +\left\Vert
(a,b,c)-(p,q,r)\right\Vert _{0},\text{ for every }(p,q,r)\in\operatorname{Gr}%
(F,G^{-1})\cap\operatorname*{cl}W. \label{ek2}%
\end{equation}

From (\ref{ek1}) we have that%
\begin{align*}
\tau(L-M^{-1})\max\{\left\Vert a-x\right\Vert ,L^{-1}\left\Vert b-y\right\Vert
,M\left\Vert c-z\right\Vert \}  &  =\left\Vert (a,b,c)-(x,y,z)\right\Vert
_{0}\\
&  \leq\left\Vert u-(y-z)\right\Vert =\left\Vert u-w\right\Vert \\
&  <\tau(L-M^{-1})\rho,
\end{align*}

\noindent hence $(a,b,c)\in W,$ and, in particular, $a\in B(x,\rho).$

If $u=b-c,$ then $u\in(F-G^{-1})(a)\subset(F-G^{-1})(B(x,\rho))$ and the
desired assertion is proved.

We want to show that $u=b-c$ is the sole possible situation. For this, suppose
by means of contradiction that $u\not =b-c.$ Fix $\varepsilon>0$ such that
$L-\varepsilon>0$ and define next%
\[
v:=(L-\varepsilon)\left\Vert u-(b-c)\right\Vert ^{-1}(u-(b-c)).
\]

\noindent Then, for every $s>0$ sufficiently small, from the $L-$ openness of
$F$ at $(a,b)\in\operatorname{Gr}F,$ we obtain that%
\[
b+sv\in B(b,Ls)\subset F(B(a,s)).
\]

Consequently, there exists $d\in B(a,s)$ such that $b+sv\in F(d).$ Obviously,
one can find $p$ with $\left\Vert p\right\Vert <1$ such that $d=a+sp.$

Also, using the $M-$ openness of $G$ at $(c,a),$ we have that for every $s>0$
sufficiently small,%
\[
B(a,s)\subset G(B(c,M^{-1}s)).
\]

Hence, one can find $e\in B(c,M^{-1}s)$ such that $d\in G(e)$ or,
equivalently, $e\in G^{-1}(d).$ Because we can write $e=c+M^{-1}sq$ with
$\left\Vert q\right\Vert <1,$ we finally have that $(a+sp,b+sv,c+M^{-1}%
sq)\in\operatorname{Gr}(F,G^{-1}).$ Taking a smaller $s$ if necessary, we also
have that $(a+sp,b+sv,c+M^{-1}sq)\in W.$ We use now (\ref{ek2}) to obtain
that, for every $s>0$ sufficiently small,%
\begin{align}
\left\Vert u-(b-c)\right\Vert  &  \leq\left\Vert u-(b+sv-c-M^{-1}%
sq)\right\Vert +s\left\Vert (p,v,M^{-1}q)\right\Vert _{0}\label{rel}\\
&  \leq\left\Vert u-(b-c)-sv\right\Vert +M^{-1}s+s\left\Vert (p,v,M^{-1}%
q)\right\Vert _{0}.\nonumber
\end{align}

But%
\[
\left\Vert u-(b-c)-sv\right\Vert =\left\vert \left\Vert u-(b-c)\right\Vert
-s(L-\varepsilon)\right\vert .
\]

Eventually for even a smaller $s,$ one obtains successively from (\ref{rel})
that%
\begin{align*}
\left\Vert u-(b-c)\right\Vert  &  \leq\left\Vert u-(b-c)\right\Vert
-s(L-\varepsilon)+M^{-1}s+s\left\Vert (p,v,M^{-1}q)\right\Vert _{0},\\
L-M^{-1}-\varepsilon &  \leq\tau(L-M^{-1})\max\{\left\Vert p\right\Vert
,L^{-1}\left\Vert v\right\Vert ,M\left\Vert M^{-1}q\right\Vert \},\\
L-M^{-1}-\varepsilon &  \leq\tau(L-M^{-1}).
\end{align*}

\noindent Passing to the limit when $\varepsilon\rightarrow0,$ we get that
$1\leq\tau,$ which is the contradiction. The proof is now complete.$\hfill
\square$

\bigskip

Let us point out that the previous theorem contains several results in
literature, especially when $F$ and $G$ are single-valued (see, for more
details, \cite[p. 412]{Urs1996}). Among these results, maybe the most famous
one is that of Graves \cite[p. 112]{Gra1950}, which can be easily deduced from
Theorem \ref{main} and Proposition \ref{link_at}. See, also, the comment after
Corollary \ref{op_fm}.

Notice that with a slight modification of the proof of Theorem \ref{main}, one
obtains the next straightforward generalization (see \cite[Theorem 3]{Urs1996}).

\begin{thm}
Let $F:X\rightrightarrows Y$ and $G_{1},...,G_{n}:Y\rightrightarrows X$ be
such that $\operatorname*{Gr}F$ and $\operatorname*{Gr}G_{i},$ $i\in
\{1,2,...,n\},$ are locally closed. Suppose $\operatorname*{Dom}(F-G_{1}%
^{-1}-...-G_{n}^{-1})$ is nonempty and let $L>0$ and $M_{1},...,M_{n}>0$ be
such that $L>M_{1}^{-1}+...+M_{n}^{-1}.$ If $F$ is $L-$open at every point
from a neighborhood of $(x,y)\in\operatorname{Gr}F$, and $G_{i}$ is $M_{i}%
-$open at every point from a neighborhood of $(z_{i},x)\in\operatorname{Gr}%
G_{i}$ for every $i\in\{1,2,...,n\}$, then $F-G_{1}^{-1}-...-G_{n}^{-1}$ is
$(L-M_{1}^{-1}-...-M_{n}^{-1})-$open at $(x,y-z_{1}-...-z_{n})$.
\end{thm}

\noindent\textbf{Proof.} Proceed as above, by taking the multifunction
$(F,G_{1}^{-1},...,G_{n}^{-1}):X\rightrightarrows Y^{n+1}$ given as
\[
(F,G_{1}^{-1},...,G_{n}^{-1})(x):=F(x)\times G_{1}^{-1}(x)\times...\times
G_{n}^{-1}(x)
\]
and observe that $(x,y,z_{1},...,z_{n})\in\operatorname*{Gr}(F,G_{1}%
^{-1},...,G_{n}^{-1}).$ The rest of the proof is more or less identical with
that of Theorem \ref{main}, using now
\[
W:=B(x,\rho)\times B(y,L\rho)\times B(z_{1},M_{1}^{-1}\rho)\times...\times
B(z_{n},M_{n}^{-1}\rho)
\]
such that $\operatorname*{Gr}(F,G_{1}^{-1},...,G_{n}^{-1})\cap
\operatorname*{cl}W$ is closed, endowing the space $X\times Y^{n+1}$ with the
norm%
\[
\left\Vert (p,q,r_{1},...,r_{n})\right\Vert _{0}:=\tau(L-M_{1}^{-1}%
-...-M_{n}^{-1})\max\{\left\Vert p\right\Vert ,L^{-1}\left\Vert q\right\Vert
,M_{1}\left\Vert r_{1}\right\Vert ,...,M_{n}\left\Vert r_{n}\right\Vert \}
\]

\noindent and applying the Ekeland variational principle to the function
$h:\operatorname{Gr}(F,G_{1}^{-1},...,G_{n}^{-1})\cap\operatorname*{cl}%
W\rightarrow\mathbb{R}_{+},$%
\begin{equation}
h(p,q,r_{1},...,r_{n}):=\left\Vert u-(q-r_{1}-...-r_{n})\right\Vert .
\tag*{$\square$}%
\end{equation}

\bigskip

The role of the next theorem is to precisely specify the constants involved in
Theorem \ref{main}. This will be a key ingredient in the proof of several
subsequent results.

\begin{thm}
\label{main_const}Let $F:X\rightrightarrows Y$ and $G:Y\rightrightarrows X$ be
two multifunctions and $(\overline{x},{\overline{y}},\overline{z})\in X\times
Y\times Y$ such that $(\overline{x},{\overline{y}})\in\operatorname{Gr}F$ and
$(\overline{z},\overline{x})\in\operatorname{Gr}G.$ Suppose that the following
assumptions are satisfied:

(i) $\operatorname{Gr}F$ is locally closed around $(\overline{x},{\overline
{y}}),$ so there exist $\alpha_{1},\beta_{1}>0$ such that $\operatorname{Gr}%
F\cap\operatorname*{cl}[B(\overline{x},\alpha_{1})\times B({\overline{y}%
},\beta_{1})]$ is closed;

(ii) $\operatorname{Gr}G$ is locally closed around $(\overline{z},\overline
{x}),$ so there exist $\alpha_{2},\beta_{2}>0$ such that $\operatorname{Gr}%
G\cap\operatorname*{cl}[B(\overline{z},\beta_{2})\times B(\overline{x}%
,\alpha_{2})]$ is closed;

(iii) there exist $L,r_{1},s_{1}>0$ such that, for every $(x^{\prime
},y^{\prime})\in\operatorname{Gr}F\cap\lbrack B(\overline{x},r_{1})\times
B({\overline{y}},s_{1})],$ $F$ is $L-$open at $(x^{\prime},y^{\prime});$

(iv) there exist $M,r_{2},s_{2}>0$ such that, for every $(v^{\prime}%
,u^{\prime})\in\operatorname{Gr}G\cap\lbrack B(\overline{z},s_{2})\times
B(\overline{x},r_{2})],$ $G$ is $M-$open at $(v^{\prime},u^{\prime});$

(v) $LM>1.$

Then for every $\rho\in]0,\varepsilon\lbrack,$ where $\varepsilon
:=\min\{\alpha_{1},\alpha_{2},L^{-1}\beta_{1},M\beta_{2},r_{1},r_{2}%
,L^{-1}s_{1},Ms_{2}\},$%
\[
B({\overline{y}}-\overline{z},(L-M^{-1})\rho)\subset(F-G^{-1})(B(\overline
{x},\rho)).
\]

Moreover, for every $\rho\in(0,2^{-1}\varepsilon)$ and every $(x,y,z)\in
\operatorname{Gr}(F,G^{-1})\cap\lbrack B(\overline{x},2^{-1}\varepsilon)\times
B({\overline{y},2}^{-1}\varepsilon)\times{B(}\overline{z},2^{-1}%
\varepsilon)],$%
\[
B(y-z,(L-M^{-1})\rho)\subset(F-G^{-1})(B(x,\rho)).
\]

\end{thm}

\noindent\textbf{Proof.} We only sketch the proof, pointing out the
differences with respect to the proof of Theorem \ref{main}.

For the first part, take $\rho\in]0,\varepsilon\lbrack,$ define, as above, the
multifunction $(F,G^{-1}),$ and observe that the choice of $\varepsilon$
implies that $\operatorname{Gr}(F,G^{-1})\cap\operatorname*{cl}W$ is closed,
where $W:=B(\overline{x},\rho)\times B({\overline{y}},L\rho)\times
B(\overline{z},M^{-1}\rho).$

Take again $u\in B({\overline{y}}-\overline{z},(L-M^{-1})\rho)$ and follow the
same steps as above to obtain that $(a,b,c)\in W.$ We only need to know that
$F$ is $L-$open at $(a,b)$ and that $G$ is $M-$open at $(c,a)$ to complete the
proof, but this follows again from the choice of $\varepsilon.$

For the second part, we define $W_{1}:=B(x,\rho)\times B(y,L\rho)\times
B(z,M^{-1}\rho)$ and we remark again that $\operatorname{Gr}(F,G^{-1}%
)\cap\operatorname*{cl}W_{1}$ is closed, because $W_{1}\subset B(\overline
{x},\min\{\alpha_{1},\alpha_{2}\})\times B({\overline{y},\beta}_{1})\times
B(\overline{z},\beta_{2}).$ The rest of the proof is the same as above,
observing only that because $(a,b,c)\in W_{1}\subset B(\overline{x}%
,\min\{r_{1},r_{2}\})\times B({\overline{y},s}_{1})\times B(\overline{z}%
,s_{2}),$ $F$ is $L-$open at $(a,b)$ and $G$ is $M-$open at $(c,a).\hfill
\square$

\bigskip

We want to emphasize that if $F$ is $L-$open around $({\overline{x}}%
,\overline{y})$, then $F$ satisfies the property from the third item, and a
similar observation is valid for $G$. Also, if one of the two multifunctions
which appear in the previous result is univoque, then one can obtain the
openness around the reference point of the difference.

\begin{cor}
\label{op_fm}Let $f:X\rightarrow Y$ be a function, $G:Y\rightrightarrows X$ be
a multifunction, $L,M>0$ and $(\overline{x},{\overline{y}},\overline{z})\in
X\times Y\times Y$ such that ${\overline{y}}=f(\overline{x})$ and
$(\overline{z},\overline{x})\in\operatorname{Gr}G.$ Suppose that the following
assumptions are satisfied:

(i) $f$ is Lipschitz continuous around $\overline{x};$

(ii) $\operatorname{Gr}G$ is locally closed around $(\overline{z},\overline
{x});$

(iii) $f$ is $L-$open around $(\overline{x},{\overline{y}});$

(iv) $G$ is $M-$open around $(\overline{z},\overline{x});$

(v) $LM>1.$

Then $f-G^{-1}$ is $(L-M^{-1})-$open around $(\overline{x},{\overline{y}%
-}\overline{z}).$
\end{cor}

\noindent\textbf{Proof.} From the Lipschitz property, we obtain the local
closedness of $\operatorname{Gr}f$ around $(\overline{x},{\overline{y}}),$ so
one deduces a similar assertion as in (i) of the previous Corollary. Observe
also that the rest of the assumptions are the same or stronger compared to the
previous result, so we suppose in the following that the conditions are
formulated using the same constants as above. Define again $\varepsilon
:=\min\{\alpha_{1},\alpha_{2},L^{-1}\beta_{1},M\beta_{2},r_{1},r_{2}%
,L^{-1}s_{1},Ms_{2}\},$ $\gamma:=\min\{4^{-1}\varepsilon,(4l)^{-1}%
\varepsilon\},$ where $l>0$ is the Lipschitz constant from (i). Take now
$\rho\in]0,2^{-1}\varepsilon\lbrack$ and $(x,v)\in\operatorname{Gr}%
(f-G^{-1})\cap\lbrack B(\overline{x},\gamma)\times B({\overline{y}%
-\overline{z},}\gamma)],$ so there exist $y=f(x)$ and $z\in G^{-1}(x)$ such
that $v=y-z.$ Then we get that%
\begin{align*}
x  &  \in B(\overline{x},\gamma)\subset B(\overline{x},2^{-1}\varepsilon),\\
y  &  \in B({\overline{y},l}\gamma)\subset B({\overline{y},4}^{-1}%
\varepsilon)\subset B({\overline{y},2}^{-1}\varepsilon),\\
z  &  \in\overline{z}+(y-{\overline{y})+}B(0,\gamma)\subset\overline
{z}+B({0,4}^{-1}\varepsilon)+B(0,4^{-1}\varepsilon)\subset B(\overline
{z},2^{-1}\varepsilon),
\end{align*}
where we used $(i)$ in the second row of inclusion.

Consequently, $(x,y,z)\in\operatorname{Gr}(f,G^{-1})\cap\lbrack B(\overline
{x},2^{-1}\varepsilon)\times B({\overline{y},2}^{-1}\varepsilon)\times
{B(}\overline{z},2^{-1}\varepsilon)],$ so using the final part of the previous
result we know that $B(y-z,(L-M^{-1})\rho)\subset(f-G^{-1})(B(x,\rho)),$ and
the proof is complete.$\hfill\square$

\bigskip

As an easy consequence of the previous corollary, one can obtain the
celebrated result of Lyusternik-Graves \cite[p. 112]{Gra1950}, taking $f$ as
the Fr\'{e}chet differential at a point of a continuously Fr\'{e}chet
differentiable function $h$ and $G:=(h-f)^{-1}.$

Another remark is that when one denotes the set-valued map $-G^{-1}$ by $F$ in
the previous corollary, then one gets (in our specific framework) the
so-called extended Lyusternik-Graves theorem (see Theorem 1 and the related
comments from \cite{DonFra2010}).

The next corollary, which can be seen as a parametric version of a result of
Graves, is the same with \cite[Proposition 3.2]{ArtMord2010} and gives
sufficient conditions for the partial metric regularity of a function, but
here is obtained as an easy consequence of Corollary \ref{op_fm}.

\begin{cor}
\label{part_open}Let $f:X\times P\rightarrow Y$ be a function between Banach
spaces which is continuous around $(\overline{x},\overline{p})$ and let
$A\in\mathcal{L}(P,Y)$ be a surjective linear operator such that there exists
$\alpha>0$ such that $\operatorname*{lop}A>\alpha$ and for every $p,p^{\prime
}\ $in a neighborhood $U$ of $\overline{p}$ and every $x\ $in a neighborhood
$V$ of $\overline{x},$%
\[
\left\Vert f(x,p)-f(x,p^{\prime})-A(p-p^{\prime})\right\Vert \leq
\alpha\left\Vert p-p^{\prime}\right\Vert .
\]

Then $f$ is open with respect to $p$ uniformly in $x$ around $((\overline
{x},\overline{p}),f(\overline{x},\overline{p}))$ with%
\[
\widehat{\operatorname*{lop}}_{p}f(\overline{x},\overline{p})\geq
\operatorname*{lop}A-\alpha.
\]

Equivalently, $f$ is metrically regular with respect to $p$ uniformly in $x$
around $((\overline{x},\overline{p}),f(\overline{x},\overline{p}))$ with%
\[
\widehat{\operatorname*{reg}}_{p}f(\overline{x},\overline{p})\leq
\frac{\operatorname*{reg}A}{1-\alpha\cdot\operatorname*{reg}A}.
\]

\end{cor}

\noindent\textbf{Proof.} Just take $x\in V$ and apply Corollary \ref{op_fm}
for $f:=-A$ and $G:=(f(x,\cdot)-A)^{-1}.$ Remark that the constants
$\varepsilon$ and $\gamma$ from the previous proof are the same for every
$x\in V$.$\hfill\square$

\bigskip

Note first that the constants in right-hand sides of the relations from the
conclusion are coming directly from the openness result, as an easy
consequence. Another remark concerns the fact that in the case where $f$ is
(strictly) partially differentiable with respect to $p,$ then, as it is often
the case in literature, one can take the partial differential with respect to
$p$ instead of $A,$ in which case $\widehat{\operatorname*{reg}}%
_{p}f(\overline{x},\overline{p})=\left\Vert (\nabla_{p}f(\overline
{x},\overline{p})^{\ast})^{-1}\right\Vert $ (see \cite{Rob1980},
\cite{ArtMord2009}, \cite[Proposition 3.4]{ArtMord2010}).

\bigskip

Here comes the second main result of the paper. For this, we use some ideas we
have previously developed in \cite{DurStr2}. Let us introduce the objects we
deal with. Remind that, for a multifunction $H:X\times P\rightrightarrows Y,$
we can define the implicit set-valued map $S:P\rightrightarrows X$ by:%
\[
S(p)=\{x\in X\mid0\in H(x,p)\}.
\]

Note that a more general solution map (see \cite{DurStr2}, \cite{NTT}) could
be investigated from the point of view of several metric regularity concepts.
We prefer the use of $S$ in the present form for clarity and unity of the results.

The result we present is an implicit multifunction theorem and shows some
interesting interrelations between the partial openness with respect to a
variable plus the Lipschitz-like property with respect to the other variable
of the original multifunction, and the Lipschitz-like, or the metric
regularity of the implicit multifunction, respectively.

\begin{thm}
\label{impl}Let $X,P$ be metric spaces, $Y$ be a normed vector space,
$H:X\times P\rightrightarrows Y$ be a set-valued map and $(\overline
{x},\overline{p},0)\in\operatorname{Gr}H$. Denote by $H_{p}(\cdot
):=H(\cdot,p),$ $H_{x}(\cdot):=H(x,\cdot).$

(i) If $H$ is open at linear rate $c>0$ with respect to $x$ uniformly in $p$
around $(\overline{x},\overline{p},0)$, then there exist $\alpha,\beta
,\gamma>0$ such that, for every $(x,p)\in B(\overline{x},\alpha)\times
B(\overline{p},\beta),$%
\begin{equation}
d(x,S(p))\leq c^{-1}d(0,H(x,p)\cap B(0,\gamma)). \label{xSpV}%
\end{equation}

If, moreover, $H$ is inner semicontinuous at $(\overline{x},\overline{p},0),$
then there exist $\alpha^{\prime},\beta^{\prime}>0$ such that, for every
$(x,p)\in B(\overline{x},\alpha^{\prime})\times B(\overline{p},\beta^{\prime
}),$%
\begin{equation}
d(x,S(p))\leq c^{-1}d(0,H(x,p)). \label{xSp}%
\end{equation}

Without the inner semicontinuity assumption on $H,$ suppose, in addition, that
$H$ is Lipschitz-like with respect to $p$ uniformly in $x$ around
$(\overline{x},\overline{p},0).$ Then $S$ is Lipschitz-like around
$(\overline{p},\overline{x})$ and%
\begin{equation}
\operatorname*{lip}S(\overline{p},\overline{x})\leq c^{-1}\widehat
{\operatorname*{lip}}_{p}H((\overline{x},\overline{p}),0). \label{lipS}%
\end{equation}

(ii) If $H$ is open at linear rate $c>0$ with respect to $p$ uniformly in $x$
around $(\overline{x},\overline{p},0)$, then there exist $\alpha,\beta
,\gamma>0$ such that, for every $(x,p)\in B(\overline{x},\alpha)\times
B(\overline{p},\beta),$%
\begin{equation}
d(p,S^{-1}(x))\leq c^{-1}d(0,H(x,p)\cap B(0,\gamma)). \label{pSxV}%
\end{equation}

Either of the following assertions are independent:

If, moreover, $H$ is inner semicontinuous at $(\overline{x},\overline{p},0),$
then there exist $\alpha^{\prime},\beta^{\prime}>0$ such that, for every
$(x,p)\in B(\overline{x},\alpha^{\prime})\times B(\overline{p},\beta^{\prime
}),$%
\begin{equation}
d(p,S^{-1}(x))\leq c^{-1}d(0,H(x,p)). \label{pSx}%
\end{equation}

Without the inner semicontinuity assumption on $H,$ suppose, in addition, that
$H$ is Lipschitz-like with respect to $x$ uniformly in $p$ around
$(\overline{x},\overline{p},0).$ Then $S$ is metrically regular around
$(\overline{p},\overline{x})$ and%
\begin{equation}
\operatorname*{reg}S(\overline{p},\overline{x})\leq c^{-1}\widehat
{\operatorname*{lip}}_{x}H((\overline{x},\overline{p}),0). \label{regS}%
\end{equation}

\end{thm}

\noindent\textbf{Proof. }We will prove only the first item, because for the
second one it suffices to observe that, defining the multifunction
$T:=S^{-1},$ the proof is completely symmetrical, using $T$ instead of $S$.
Moreover, using Proposition \ref{link_around}, we know that
$\operatorname*{reg}S(\overline{p},\overline{x})=\operatorname*{lip}%
T(\overline{x},\overline{p})$ and then (\ref{regS}) follows from (\ref{lipS}).

For the (i) item, we know that there exist $r,s,t,c,\varepsilon>0$ such that,
for every $\rho\in(0,\varepsilon),$ every $p\in B(\overline{p},t)$ and every
$(x,y)\in\operatorname{Gr}H_{p}\cap\lbrack B(\overline{x},r)\times B(0,s)],$%
\[
B(y,c\rho)\subset H_{p}(B(x,\rho)).
\]

Take now $\rho\in]0,\min\{\varepsilon,c^{-1}s\}[.$ Set $\alpha:=r,$
$\beta:=t,$ $\gamma:=c\rho$ and fix arbitrary $(x,p)\in B(\overline{x}%
,\alpha)\times B(\overline{p},\beta).$ If $H(x,p)\cap B(0,c\rho)=\emptyset,$
then $d(0,H(x,p)\cap B(0,\gamma))=+\infty$ and (\ref{xSpV}) trivially holds.
Suppose next that $H(x,p)\cap B(0,c\rho)\not =\emptyset.$ If $0\in H(x,p),$
then $0\in H(x,p)\cap B(0,c\rho),$ and, again, (\ref{xSpV}) trivially holds.
Suppose now that $0\not \in H(x,p)\cap B(0,c\rho).$ Then for every $\xi>0,$
there exists $y_{\xi}\in H(x,p)\cap B(0,c\rho)$ such that%
\[
\left\Vert y_{\xi}\right\Vert <d(0,H(x,p)\cap B(0,c\rho))+\xi.
\]

Because $d(0,H(x,p)\cap B(0,c\rho))<c\rho,$ we can choose $\xi$ sufficiently
small such that $d(0,H(x,p)\cap B(0,c\rho))+\xi<c\rho.$ Consequently,%
\begin{equation}
0\in B(y_{\xi},d(0,H(x,p)\cap B(0,c\rho))+\xi)\subset B(y_{\xi},c\rho).
\label{0inB}%
\end{equation}

Observe now that $x\in B(\overline{x},r),$ $p\in B(\overline{p},t),$ $y_{\xi
}\in B(0,d(0,H(x,p)\cap B(0,c\rho))+\xi)\subset B(0,c\rho)\subset B(0,s),$
$y_{\xi}\in H(x,p)$ and denote $\rho_{0}:=c^{-1}(d(0,H(x,p)\cap B(0,c\rho
))+\xi)<\rho<\varepsilon.$

But we know that%
\[
B(y_{\xi},c\rho_{0})\subset H_{p}(B(x,\rho_{0})),
\]

\noindent hence, using also (\ref{0inB}), one obtains that there exists
$x_{0}\in B(x,\rho_{0})$ such that $0\in H(x_{0},p),$ which is equivalent to
$x_{0}\in S(p).$ Then%
\[
d(x,S(p))\leq d(x,x_{0})<\rho_{0}=c^{-1}(d(0,H(x,p)\cap B(0,c\rho))+\xi).
\]

Making $\xi\rightarrow0,$ we obtain (\ref{xSpV}).

In the case that $H$ is inner semicontinuous at $(\overline{x},\overline
{p},0),$ one can find $\delta_{1},\delta_{2}>0$ such that, for every $(x,p)\in
B(\overline{x},\delta_{1})\times B(\overline{p},\delta_{2}),$%
\begin{equation}
H(x,p)\cap B(0,c\rho)\not =\emptyset. \label{iscH}%
\end{equation}

Then $\alpha^{\prime}:=\min\{r,\gamma\},$ $\beta^{\prime}:=\min\{t,\delta\}$
are appropriate such that (\ref{xSp}) is satisfied.

Suppose now that $H$ is Lipschitz-like with respect to $p$ uniformly in $x$
around $(\overline{x},\overline{p},0).$ Then there exist $l,a,b,\tau>0$ such
that $\tau<c\rho$ and for every $x\in B(\overline{x},a)$ and every
$p_{1},p_{2}\in B(\overline{p},b),$%
\begin{equation}
H(x,p_{1})\cap D(0,\tau)\subset H(x,p_{2})+ld(p_{1},p_{2})\mathbb{D}_{Y}.
\label{lipHp}%
\end{equation}

Take $\overline{\alpha}:=\min\{a,\alpha\},$ $\overline{\beta}=\min
\{b,\beta,(2l)^{-1}\tau\},$ $p_{1},p_{2}\in B(\overline{p},\overline{\beta})$
and $x\in S(p_{1})\cap B(\overline{x},\overline{\alpha}).$ Then $0\in
H(x,p_{1})\cap D(0,\tau),$ whence, using (\ref{lipHp}), there exists
$y^{\prime\prime}\in\mathbb{D}_{Y}$ such that $y^{\prime}:=l\cdot
d(p_{1},p_{2})y^{\prime\prime}\in H(x,p_{2})$ with $\left\Vert y^{\prime
}\right\Vert \leq l\cdot\lbrack d(p_{1},\overline{p})+d(\overline{p}%
,p_{2})]\leq\tau<c\rho.$ Hence, $y^{\prime}\in H(x,p_{2})\cap B(0,\gamma),$ so
using (\ref{xSpV}), we get that%
\[
d(x,S(p_{2}))\leq c^{-1}d(0,H(x,p_{2})\cap B(0,\gamma))\leq c^{-1}\left\Vert
y^{\prime}\right\Vert \leq c^{-1}ld(p_{1},p_{2}).
\]

Consequently, because $l$ can be chosen arbitrarily close to $\widehat
{\operatorname*{lip}}_{p}H((\overline{x},\overline{p}),0)),$ it follows that
$S$ is Lipschitz-like around $(\overline{p},\overline{x})$ and
$\operatorname*{lip}S(\overline{p},\overline{x})\leq c^{-1}\widehat
{\operatorname*{lip}}_{p}H((\overline{x},\overline{p}),0)).$ The proof is now
complete.$\hfill\square$

\bigskip

It is worth to be mentioned that one can obtain in the final parts of the two
items of the above theorem even a kind of graphical regularity, following the
technique from \cite[Theorem 5.2]{DurStr2} (see, also, \cite{DurStr3}).

\section{Applications}

This section is dedicated to the investigation of the case where the mapping
$H$ is given as a sum of two set-valued maps in the sense we shall precise.
However, we start with an application of Theorem \ref{impl} and we get an
implicit multifunction result which generalizes \cite[Theorem 3.5]%
{ArtMord2010} to the case where the set-valued map $\Gamma$ is constructed
using the sum between a function and a multifunction. Also, in the virtue of
Corollary \ref{part_open}, we can conclude that this result generalizes also
\cite[Lemma 3.1]{ArtMord2009}.

\begin{pr}
\label{gama}Let $X,Y,Z,W$ be Banach spaces, $F:X\times Y\rightrightarrows Z$
be a multifunction, $g:W\rightarrow Z$ be a function and $(\overline
{x},{\overline{y},}\overline{z},\overline{w})\in X\times Y\times Z\times W$ be
such that $\overline{z}:=-g(\overline{w})\in F(\overline{x},{\overline{y}}).$
Consider the implicit multifunction $\Gamma:Y\times W\rightrightarrows X$
defined by%
\[
\Gamma(y,w):=\{x\in X\mid0\in F(x,y)+g(w)\}.
\]

Suppose that the following conditions are satisfied:

(i) $F$ is Lipschitz-like with respect to $y$ uniformly in $x$ around
$((\overline{x},{\overline{y}),}\overline{z})$ with constant $\eta\geq0;$

(ii) $F$ is metrically regular with respect to $x$ uniformly in $y$ around
$((\overline{x},{\overline{y}),}\overline{z})$ with constant $k>0;$

(iii) $g$ is locally Lipschitzian around $\overline{w}$ with constant
$\lambda>0.$

Then there exists $\alpha>0$ such that for every $(y,w),(y^{\prime},w^{\prime
})\in D({\overline{y},}\alpha)\times D(\overline{w},\alpha)$ and for every
$\varepsilon>0,$%
\[
\Gamma(y^{\prime},w^{\prime})\cap D(\overline{x},\alpha)\subset\Gamma
(y,w)+(k+\varepsilon)(\eta\left\Vert y-y^{\prime}\right\Vert +\lambda
\left\Vert w-w^{\prime}\right\Vert )\mathbb{D}_{X}.
\]

In particular, $\Gamma$ is Lipschitz-like around $(({\overline{y},}%
\overline{w}),\overline{x})$ with the following estimate%
\[
\operatorname*{lip}\Gamma(({\overline{y},}\overline{w}),\overline{x}%
)\leq\widehat{\operatorname*{reg}}_{x}F((\overline{x},{\overline{y}%
),}\overline{z})\cdot\max\{\widehat{\operatorname*{lip}}_{y}F(\overline
{x},{\overline{y},}\overline{z}),\operatorname*{lip}g(\overline{w})\}.
\]

\end{pr}

\noindent\textbf{Proof.} Define $P:=Y\times W$ and $H:X\times
P\rightrightarrows Z$ by $H(x,(y,w)):=F(x,y)+g(w).$ Observe also that%
\[
\Gamma(y,w)=\{x\in X\mid0\in H(x,(y,w))\},
\]
and denote $H_{(y,w)}(\cdot):=H(\cdot,(y,w)).$ For every $y$ close to
${\overline{y}}$ we know from (ii) and Proposition \ref{link_around} that
$F_{y}$ is $k^{-1}-$open at points from its graph around $(\overline
{x},\overline{z}).$ Therefore, we can conclude that there exists
$\varepsilon>0$ such that for every $\rho\in]0,\varepsilon\lbrack,$ every
$y\in B({\overline{y},}\varepsilon)$ and every $(x,z^{\prime})\in
\operatorname{Gr}F_{y}\cap\lbrack B(\overline{x},\varepsilon)\times
B(\overline{z},\varepsilon)],$
\begin{equation}
B(z^{\prime},k^{-1}\rho)\subset F_{y}(B(x,\rho)). \label{Fy}%
\end{equation}

Take $(y,w)\in B({\overline{y},}\varepsilon)\times B({\overline{w},2}%
^{-1}\lambda^{-1}\varepsilon)$ and $(x,z)\in\operatorname{Gr}H_{(y,w)}%
\cap\lbrack B(\overline{x},\varepsilon)\times B({0,2}^{-1}\varepsilon)],$ and
$z^{\prime}:=z-g(w).$ Then, using (iii), we have that%
\[
\left\Vert z^{\prime}-\overline{z}\right\Vert \leq\left\Vert z\right\Vert
+\left\Vert \overline{z}+g(w)\right\Vert <2^{-1}\varepsilon+\lambda\left\Vert
w-\overline{w}\right\Vert <\varepsilon,
\]
so one can use (\ref{Fy}) to prove that%
\[
B(z,k^{-1}\rho)=B(z^{\prime},k^{-1}\rho)-g(w)\subset F_{y}(B(x,\rho
))-g(w)=H_{(y,w)}(B(x,\rho)).
\]
But this shows, applying Theorem \ref{impl}, that there exist $\beta>0$ such
that, for every $(y,z,w)\in B(\overline{y},\beta)\times B(\overline{z}%
,\beta)\times B(\overline{w},\beta),$%
\begin{equation}
d(y,\Gamma(y,w))\leq kd(0,H(x,(y,w))\cap B(0,\beta)). \label{relGa}%
\end{equation}

We want to prove that there exists $c>0$ such that $(\eta+1)c<\beta$, for
every $x\in B(\overline{x},c),$ $(y,w),(y^{\prime},w^{\prime})\in
B({\overline{y},}c)\times B(\overline{w},c),$%

\begin{equation}
H(x,(y,w))\cap D(0,c)\subset H(x,(y^{\prime},w^{\prime}))+(\eta\left\Vert
y-y^{\prime}\right\Vert +\lambda\left\Vert w-w^{\prime}\right\Vert
)\mathbb{D}_{Z}. \label{relH}%
\end{equation}

In particular, we will prove that $H$ is Lipschitz-like with respect to
$(y,w)$ uniformly in $x$ around $(\overline{x},({\overline{y},}\overline
{z}),0).$

Because of (i), we know that there exists $a>0$ such that for every $x\in
B(\overline{x},a)$ and every $y,y^{\prime}\in B({\overline{y},a}),$%
\begin{equation}
F(x,y)\cap D(\overline{z},a)\subset F(x,y^{\prime})+\eta\left\Vert
y-y^{\prime}\right\Vert \mathbb{D}_{Z}. \label{relF}%
\end{equation}

Also, because of (iii), we can find $b>0$ such that for every $w,w^{\prime}\in
B(\overline{w},b),$%
\begin{equation}
g(w)\in g(w^{\prime})+\lambda\left\Vert w-w^{\prime}\right\Vert \mathbb{D}%
_{Z}. \label{lipg}%
\end{equation}

Choose now $c\in]0,\min\{(\lambda+1)^{-1}a,b,(\eta+\lambda)^{-1}\beta\}[$ and
take arbitrary $x\in B(\overline{x},c),$ $(y,w),(y^{\prime},w^{\prime})\in
B({\overline{y},}c)\times B(\overline{w},c).$ Furthermore, choose $z\in
H(x,(y,w))\cap D(0,c).$ Then $z-g(w)\in F(x,y)$ and because of (\ref{lipg}),
we know that $-g(w)\in D(\overline{z},\lambda\left\Vert w-\overline
{w}\right\Vert )\subset B(\overline{z},\lambda c),$ whence $z-g(w)\in
B(\overline{z},(\lambda+1)c)\subset B(\overline{z},a).$ One can use now
(\ref{relF}) to obtain that $z-g(w)\in F(x,y^{\prime})+\eta\left\Vert
y-y^{\prime}\right\Vert \mathbb{D}_{Z}.$ Adding with (\ref{lipg}), one finally
gets (\ref{relH}).

Take now $\gamma\in]0,\min\{\beta,2^{-1}c\}[,$ $(y,w),(y^{\prime},w^{\prime
})\in B({\overline{y},}\gamma)\times B(\overline{w},\gamma),$ $x\in
\Gamma(y,w)\cap D(\overline{x},\gamma)$ and arbitrary $\varepsilon>0.$ Hence,
$0\in H(x,(y,w))\cap D(0,c).$ Then, using (\ref{relH}), we have that there
exists $u\in\mathbb{D}_{Z}$ such that $(\eta\left\Vert y-y^{\prime}\right\Vert
+\lambda\left\Vert w-w^{\prime}\right\Vert )u\in H(x,(y^{\prime},w^{\prime
}))\cap B(0,\beta)$ (because $(\eta+\lambda)c<\beta$). Finally, using
(\ref{relGa}), we get that
\begin{align*}
d(x,\Gamma(y^{\prime},w^{\prime}))  &  \leq kd(0,H(x,(y^{\prime},w^{\prime
}))\cap B(0,\beta))\leq k(\eta\left\Vert y-y^{\prime}\right\Vert
+\lambda\left\Vert w-w^{\prime}\right\Vert )\\
&  <(k+\varepsilon)(\eta\left\Vert y-y^{\prime}\right\Vert +\lambda\left\Vert
w-w^{\prime}\right\Vert ),
\end{align*}

\noindent which completes the proof.$\hfill\square$

\bigskip

Note that the final conclusion of Theorem \ref{impl} (i), but the estimations
(\ref{xSpV}) and (\ref{xSp}), could be obtained as a consequence of
Proposition \ref{gama}, taking $g\equiv0$. We note as well that for $W:=Z$ and
$g(z):=-z$ for every $z$\ we can get an even more general implicit
multifunction result (see \cite[Theorem 5.2]{DurStr2}).

\bigskip

The next technical notion will be used in the sequel, mainly to prove a
Lipschitz-like property of the sum between two multifunctions. In this way (in
contrast to \cite[Corollary 18]{NTT}), we avoid the strong requirements of the
single-valuedness and full Lipschitz property of the field map $G$ at the
reference point.

\begin{df}
\label{sum-sta}Let $F:X\rightrightarrows Y,$ $G:X\rightrightarrows Y$ be two
multifunctions and $(\overline{x},{\overline{y},}\overline{z})\in X\times
Y\times Y$ such that ${\overline{y}\in F(\overline{x}),}$ ${\overline{z}\in
G(\overline{x}).}$ We say that the multifunction $(F,G)$ is locally sum-stable
around $(\overline{x},{\overline{y},}\overline{z})$ if for every
$\varepsilon>0$ there exists $\delta>0$ such that, for every $x\in
B(\overline{x},\delta)$ and every $w\in(F+G)(x)\cap B({\overline{y}+}%
\overline{z},\delta),$ there exist $y\in F(x)\cap B({\overline{y},}%
\varepsilon)$ and $z\in G(x)\cap B(\overline{z},\varepsilon)$ such that
$w=y+z.$
\end{df}

This definition is illustrated, at a first glance, by two simple examples.
First example displays a simple situation where this condition holds true.
Note that $\mathbb{R}$ and $\mathbb{Q}$ denote the fields of reals and
rationals, respectively.

\begin{examp}
Let $F:\mathbb{R\rightrightarrows R}$ given, for any $x\in\mathbb{R},$ by
\[
F(x):=]0,\left\vert x\right\vert [.
\]
Take now $G:=F.$ We will prove that $(F,G)$ is locally sum-stable around any
$(\overline{x},{\overline{y}},\overline{z})\in\operatorname*{Gr}(F,G).$ In
order to prove this, let us take, without loosing of generality, $\overline
{x}>0$ and $\overline{y},\overline{z}\in]0,\overline{x}[.$ Then fix
$\varepsilon>0$ and choose $\delta\in]0,\min\{2^{-1}(\overline{x}-\overline
{y}),2^{-1}(\overline{x}-\overline{z}),\varepsilon,{\overline{y},}\overline
{z}\}[.$ Further, choose $x\in]\overline{x}-\delta,\overline{x}+\delta
\lbrack\ $and $w\in(F+G)(x)\cap B({\overline{y}+}\overline{z},\delta
)=]0,2x[\cap]{\overline{y}+}\overline{z}-\delta,{\overline{y}+}\overline
{z}+\delta\lbrack.$ Now, if $w<{\overline{y}+}\overline{z},$ then there exist
$\gamma\in]0,\delta\lbrack$ s.t. $w={\overline{y}+}\overline{z}-\gamma.$
Consider $y:=\overline{y}-2^{-1}\gamma$ and $z=\overline{z}-2^{-1}\gamma.$
Then, $\overline{y}-y=2^{-1}\gamma<\delta<\varepsilon,$ and $y>\overline
{y}-\delta>0.$ Moreover, $y=\overline{y}-2^{-1}\gamma<\overline{x}-2\delta<x,$
i.e. $y\in F(x).$ Similarly for $z.$ If $w>{\overline{y}+}\overline{z},$ then
there exist $\gamma\in]0,\delta\lbrack$ s.t.  $w={\overline{y}+}\overline
{z}+\gamma.$ Consider $y:=\overline{y}+2^{-1}\gamma$ and $z=\overline
{z}+2^{-1}\gamma.$ Of course, $y>0$ and $y-\overline{y}=2^{-1}\gamma
<\delta<\varepsilon$. On the other hand, $y=\overline{y}+2^{-1}\gamma
<\overline{y}+\delta<\overline{x}-\delta<x,$ i.e. $y\in F(x).$ The same
calculation holds for $z$ and the proof is complete.
\end{examp}

The next example describes a situation where the sum-stable condition does not hold.

\begin{examp}
Let $F:\mathbb{R\rightrightarrows R}$ be given by
\[
F(x):=\left\{
\begin{array}
[c]{ll}%
\lbrack-1,1], & \text{if }x\in\mathbb{R}\setminus\{1\}\\
2, & \text{if }x=1
\end{array}
\right.
\]
and $G:\mathbb{R\rightrightarrows R}$ given, for any $x\in\mathbb{R},$ by
$G(x):=-F(x).$ Then one can easily see that $(F,G)$ is not sum-stable at
$(\overline{x},{\overline{y},}\overline{z}):=(1,2,-2)$ because $0\in(F+G)(x)$
for every $x\in\mathbb{R}$ but if $x\neq1,$ then we cannot write $0$ as a sum
between an element in $F(x)\cap B(2,2^{-1})$ and an element in $G(x)\cap
B(-2,2^{-1}).$
\end{examp}

Remind that a\ multifunction $F:X\rightrightarrows Y$ is said to be Lipschitz
around $\overline{x}\in\operatorname*{Dom}F$ with constant $L>0$ if there
exists a neighborhood $U\in\mathcal{V}(\overline{x})$ such that, for every
$x,u\in U,$%
\begin{equation}
F(x)\subset F(u)+L\left\Vert x-u\right\Vert \mathbb{D}_{Y}. \label{lip}%
\end{equation}

Of course, this property is (much) stronger than the Lipschitz-like property,
having the great advantage to be stable at summation. More precisely, if $F,G$
are two multifunctions which are Lipschitz around some point $\overline{x}%
\in\operatorname*{Dom}F\cap\operatorname*{Dom}G,$ then $F+G$ is Lipschitz
around $\overline{x}.$

The next, more elaborated example, shows that the Lipschitz property of both
multifunctions does not ensure the sum-stable property.

\begin{examp}
Let $F:\mathbb{R\rightrightarrows R}$ be given by $F(x):=[1,+\infty\lbrack
\cap\mathbb{Q}$ for every $x\in\mathbb{R}$ and $G:\mathbb{R\rightrightarrows
R}$ given by
\[
G(x):=\left\{
\begin{array}
[c]{ll}%
]-\infty,-1], & \text{if }x\in\mathbb{R}\setminus\{1-\frac{1}{n}\mid
n\in\mathbb{N}\setminus\{0\}\}\\
]-\infty,-1]\cap\mathbb{Q\cup\{-}2\mathbb{+}\frac{\sqrt{2}}{n}\mathbb{\}}%
\text{,} & \text{if }x=1-\frac{1}{n},\text{ }n\in\mathbb{N}\setminus\{0\}.
\end{array}
\right.
\]
It is easy to verify that both $F$ and $G$ are Lipschitz around $\overline
{x}:=1,$ because, for example, if we add to the set $]-\infty,-1]\cap
\mathbb{Q\cup\{-}2\mathbb{+}\frac{\sqrt{2}}{n}\mathbb{\}}$ any ball, then we
cover all the interval $]-\infty,-1].$ Take now $\overline{y}:=1$ and
$\overline{z}:=-1,$ fix $\varepsilon:=2^{-1}$ and pick arbitrary $\delta>0$.
Then choose $x_{n}:=1-\frac{1}{n},$ $n\in\mathbb{N}\setminus\{0\}$ and observe
that $x_{n}\in]1-\delta,1+\delta\lbrack$ and
\[
w_{n}:=\frac{\sqrt{2}}{n}\in(F+G)(x_{n})\cap B(0,\delta)
\]
for any $n$ sufficiently large$.$ But $w_{n}$ can be obtained only by the sum
between $2\in F(x_{n})$ and $\mathbb{-}2\mathbb{+}\frac{\sqrt{2}}{n}\in
G(x_{n})$ and both these values are not in the balls $B(\overline{y},2^{-1})$
and $B(\overline{z},2^{-1})$ respectively. Therefore, $(F,G)$ is not
sum-stable at $(\overline{x},\overline{y},\overline{z})=(1,1,-1).$
\end{examp}

Next proposition indicates a first general situation where the local-sum
stability holds. Recall that a function $f:X\rightarrow Y$ is calm at
$\overline{x}$ if there exist $\alpha,l>0$ such that, for every $x\in
B(\overline{x},\alpha),$%
\[
\left\Vert f(x)-f(\overline{x})\right\Vert \leq l\left\Vert x-\overline
{x}\right\Vert .
\]

\begin{pr}
\label{fGsta}Let $f:X\rightarrow Y$ be a function, $G:X\rightrightarrows Y$ be
a multifunction and $\overline{x}\in X$ such that $0\in f(\overline
{x})+G(\overline{x}).$ If $f$ is calm at $\overline{x},$ then $(f,G)$ is
locally sum-stable around $(\overline{x},f({\overline{x}),-f(\overline{x})}).$
\end{pr}

\noindent\textbf{Proof. }Suppose that the constants from the calmness property
of $f$ are the same as above, take arbitrarily $\varepsilon>0$ and choose
$\delta\in]0,\min\{2^{-1}\varepsilon,(2l)^{-1}\varepsilon,\alpha\}[.$ Pick now
$x\in B(\overline{x},\delta)$ and $w\in(f+G)(x)\cap B({0},\delta).$ Then
$w-f(x)\in G(x).$ Moreover, $\left\Vert f(x)-f(\overline{x})\right\Vert \leq
l\left\Vert x-\overline{x}\right\Vert <l\delta<2^{-1}\varepsilon,$ so $f(x)\in
B(f(\overline{x}),\varepsilon).$ Also, $w-f(x)\in B({0},\delta)+B(-f(\overline
{x}),2^{-1}\varepsilon)\subset B(-f(\overline{x}),\varepsilon).$ The proof is
now complete.$\hfill\square$

\bigskip

The next lemma is the main motivation for introducing the sum-stable property.

\begin{lm}
\label{lip_H}Let $F:X\rightrightarrows Y,$ $G:X\rightrightarrows Y$ be two
multifunctions. Suppose that $F$ is Lipschitz-like around $(\overline
{x},{\overline{y}})\in\operatorname{Gr}F,$ that $G$ is Lipschitz-like around
$(\overline{x},{\overline{z}})\in\operatorname{Gr}G$ and that $(F,G)$ is
locally sum-stable around $(\overline{x},{\overline{y},}\overline{z}).$ Then
the multifunction $F+G$ is Lipschitz-like around $(\overline{x},{\overline
{y}+}\overline{z}).$ Moreover, the following relation holds true%
\begin{equation}
\operatorname*{lip}(F+G)(\overline{x},{\overline{y}+}\overline{z}%
)\leq\operatorname*{lip}F(\overline{x},{\overline{y}})+\operatorname*{lip}%
G(\overline{x},\overline{z}). \label{lip_sum}%
\end{equation}

\end{lm}

\noindent\textbf{Proof. }Using the Lipschitz-like properties of $F$ and $G,$
one can find $\alpha,l,k>0$ such that, for every $x,u\in B(\overline{x}%
,\alpha),$%
\begin{align}
F(x)\cap B({\overline{y},}\alpha)  &  \subset F(u)+l\left\Vert x-u\right\Vert
\mathbb{D}_{Y},\label{lipF}\\
G(x)\cap B({\overline{z},}\alpha)  &  \subset G(u)+k\left\Vert x-u\right\Vert
\mathbb{D}_{Y}. \label{lipG}%
\end{align}
But using the local sum-stability for $\varepsilon:=\alpha>0,$ we can find
$\delta\in]0,\alpha\lbrack$ such that, for every $x\in B(\overline{x},\delta)$
and every $w\in(F+G)(x)\cap B({\overline{y}+}\overline{z},\delta),$ there
exist $y\in F(x)\cap B({\overline{y},}\alpha)$ and $z\in G(x)\cap
B(\overline{z},\alpha)$ such that $w=y+z.$ Consequently, using (\ref{lipF})
and (\ref{lipG}), for every $u\in B(\overline{x},\delta),$ $w\in
(F+G)(u)+(l+k)\left\Vert x-u\right\Vert \mathbb{D}_{Y}.$ The relation
(\ref{lip_sum}) follows from the fact that constants $l$ and $k$ can be chosen
arbitrarily close to $\operatorname*{lip}F(\overline{x},{\overline{y}})$ and
$\operatorname*{lip}G(\overline{x},\overline{z}),$ respectively.$\hfill
\square$

\bigskip

We would like to continue with our examples above in order to illustrate the
fact that sum-stable property is essential in Lemma \ref{lip_H}. Basically, we
need an example of two Lipschitz-like multifunctions for which the sum is not
Lipschitz-like at the reference point.

\begin{examp}
Let $F:\mathbb{R\rightrightarrows R}$ be given by
\[
F(x):=\left\{
\begin{array}
[c]{ll}%
\lbrack1,2], & \text{if }x\in\mathbb{R}\setminus\{1\}\\
\lbrack1,2]\cup\{0\}, & \text{if }x=1
\end{array}
\right.
\]
and $G:\mathbb{R\rightrightarrows R}$ given by $G(x):=[1,2]$ for every
$x\in\mathbb{R}$.

As above, it is easy to see that both $F$ and $G$ are Lipschitz-like around
$(\overline{x},\overline{y})=(1,1)$ and $(\overline{x},\overline{z})=(1,1),$
respectively. But the multifunction $F+G:\mathbb{R}\rightrightarrows
\mathbb{R}$, given by%
\[
(F+G)(x)=\left\{
\begin{array}
[c]{ll}%
\lbrack2,4], & \text{if }x\in\mathbb{R}\setminus\{1\}\\
\lbrack1,4], & \text{if }x=1
\end{array}
\right.  ,
\]

\noindent is not Lipschitz-like around $(\overline{x},\overline{y}%
+\overline{z})=(1,2).$ Indeed, suppose by contradiction that there exists
$L>0$ and $\alpha\in]0,\min\{1,L^{-1}\}[$ such that for any $x,u\in
\lbrack1-\alpha,1+\alpha]$%
\begin{equation}
(F+G)(x)\cap\lbrack2-\alpha,2+\alpha]\subset(F+G)(u)+L\left\vert
x-u\right\vert [-1,1]. \label{rpL}%
\end{equation}
Consider now $x:=1$ and $u:=1-\alpha^{2}$ such that $x,u\in\lbrack
1-\alpha,1+\alpha].$ Clearly,%
\[
2-\alpha\in(F+G)(x)\cap\lbrack2-\alpha,2+\alpha].
\]
Following (\ref{rpL}), we should have:%
\[
2-\alpha\in(F+G)(1-\alpha^{2})+L\alpha^{2}[-1,1]
\]
and, in particular,
\[
2-\alpha\in\lbrack2-L\alpha^{2},4+L\alpha^{2}].
\]
But this requires that%
\[
\alpha\leq L\alpha^{2},
\]
which contradicts the choice of $\alpha$. The contradiction shows that we
cannot have the Lipschitz-like property of the sum.

Observe that, in the virtue of Lemma \ref{lip_H}, $(F,G)$ cannot be
locally-sum stable around $(1,1,1).$ Indeed, take $\varepsilon\in]0,2^{-1}[.$
Then for every $\delta>0,$ choose $n\in\mathbb{N}$ such that $n>\max
\{\delta,1\}.$ Taking now $x_{\delta}:=1\in]1-\delta,1+\delta\lbrack$ and
$w_{\delta}:=2-n^{-1}\delta\in(F+G)(x_{\delta})\cap]2-\delta,2+\delta\lbrack,$
one can easily see that, for every $y\in F(x_{\delta})\cap]1-\varepsilon
,1+\varepsilon\lbrack=[1,1+\varepsilon\lbrack$ and every $z\in G(x_{\delta
})\cap]1-\varepsilon,1+\varepsilon\lbrack=[1,1+\varepsilon\lbrack,$
$w_{\delta}<y+z.$
\end{examp}

Next, we adapt the definition of local-sum stability to the parametric case,
in order to use this notion in the general context of variational systems.

\begin{df}
\label{sum-sta_param}Let $F:X\times P\rightrightarrows Y,$
$G:X\rightrightarrows Y$ be two multifunctions and $(\overline{x},\overline
{p},{\overline{y},}\overline{z})\in X\times P\times Y\times Y$ such that
${\overline{y}\in F(\overline{x},\overline{p}),}$ ${\overline{z}\in
G(\overline{x}).}$ We say that the multifunction $(F,G)$ is locally sum-stable
around $(\overline{x},\overline{p},{\overline{y},}\overline{z})$ if for every
$\varepsilon>0$ there exists $\delta>0$ such that, for every $(x,p)\in
B(\overline{x},\delta)\times B(\overline{p},\delta)$ and every $w\in
(F_{p}+G)(x)\cap B({\overline{y}+}\overline{z},\delta),$ there exist $y\in
F_{p}(x)\cap B({\overline{y},}\varepsilon)$ and $z\in G(x)\cap B(\overline
{z},\varepsilon)$ such that $w=y+z.$
\end{df}

Similarly to Proposition \ref{fGsta}, one can easily prove the next (adapted) result.

\begin{pr}
Let $f:X\times P\rightarrow Y$ be a function, $G:X\rightrightarrows Y$ be a
multifunction and $(\overline{x},\overline{p})\in X\times P$ such that $0\in
f(\overline{x},\overline{p})+G(\overline{x}).$ If $f$ is calm at
$(\overline{x},\overline{p}),$ then $(f,G)$ is locally sum-stable around
$(\overline{x},\overline{p},f({\overline{x},\overline{p}),-f(\overline
{x},\overline{p})}).$
\end{pr}

Also, Lemma \ref{lip_H} has the following variant in the parametric case.

\begin{lm}
\label{lip_H_par}Let $F:X\times P\rightrightarrows Y,$ $G:X\rightrightarrows
Y$ be two multifunctions. Suppose that $F$ is Lipschitz-like with respect to
$x$ uniformly in $p$ around $((\overline{x},\overline{p}),{\overline{y}}%
)\in\operatorname{Gr}F,$ that $G$ is Lipschitz-like around $(\overline
{x},{\overline{z}})\in\operatorname{Gr}G$ and that $(F,G)$ is locally
sum-stable around $(\overline{x},\overline{p},{\overline{y},}\overline{z}).$
Then the multifunction $H:X\times P\rightrightarrows Y$ given by
$H(x,p):=F(x,p)+G(x)$ is Lipschitz-like with respect to $x$ uniformly in $p$
around $((\overline{x},\overline{p}),{\overline{y}+}\overline{z}).$ Moreover,
the following relation holds true%
\begin{equation}
\widehat{\operatorname*{lip}}_{x}H((\overline{x},\overline{p}),{\overline{y}%
+}\overline{z})\leq\widehat{\operatorname*{lip}}_{x}F((\overline{x}%
,\overline{p}),{\overline{y}})+\operatorname*{lip}G(\overline{x},\overline
{z}). \label{lip_sum_par}%
\end{equation}

\end{lm}

\bigskip

The following result deduces the metric regularity of $S$ under appropriate
assumptions on the multifunctions $F$ and $G,$ which arrive naturally from
Theorem \ref{impl}. Namely, part of these assumptions are stated in order to
ensure the Lipschitz-like property of the sum multifunction $H$ with respect
to $x$ uniformly in $p$ around $((\overline{x},\overline{p}),{0}),$ with
(\ref{lip_sum_par}) satisfied. Using (\ref{regS}), we expect to have that
$\operatorname*{reg}S(\overline{p},\overline{x})\leq c^{-1}\widehat
{\operatorname*{lip}}_{x}H((\overline{x},\overline{p}),0),$ where $c>0$ is the
rate of linear openness with respect to $p$ of $H,$ but is easy to see that
$c$ must be $\widehat{\operatorname*{lop}}_{p}F((\overline{x},\overline
{p}),\overline{y})=$ $\left(  \widehat{\operatorname*{reg}}_{p}F((\overline
{x},\overline{p}),\overline{y})\right)  ^{-1}.$ Therefore, (\ref{rS}) below
holds in a natural way.

\begin{thm}
\label{mreg_sol}Let $X,Y,P$ be Banach spaces, $F:X\times P\rightrightarrows
Y,$ $G:X\rightrightarrows Y$ be two set-valued maps and $(\overline
{x},\overline{p},\overline{y})\in X\times P\times Y$ such that $\overline
{y}\in F(\overline{x},\overline{p})$ and $-\overline{y}\in G(\overline{x})$.
Suppose that the following assumptions are satisfied:

(i) $(F,G)$ is locally sum-stable around $(\overline{x},\overline
{p},{\overline{y},}-{\overline{y}});$

(ii) $F$ is Lipschitz-like with respect to $x$ uniformly in $p$ around
$((\overline{x},\overline{p}),\overline{y});$

(iii) $F$ is metrically regular with respect to $p$ uniformly in $x$ around
$((\overline{x},\overline{p}),\overline{y});$

(iv) $G$ is Lipschitz-like around $(\overline{x},-{\overline{y}}).$

Then $S$ is metrically regular around $(\overline{p},\overline{x}).$ Moreover,
the next relation holds%
\begin{equation}
\operatorname*{reg}S(\overline{p},\overline{x})\leq\widehat
{\operatorname*{reg}}_{p}F((\overline{x},\overline{p}),\overline{y}%
)\cdot\lbrack\widehat{\operatorname*{lip}}_{x}F((\overline{x},\overline
{p}),\overline{y})+\operatorname*{lip}G(\overline{x},-{\overline{y}})].
\label{rS}%
\end{equation}

\end{thm}

\noindent\textbf{Proof. }Define $H:X\times P\rightrightarrows Y$ by
\begin{equation}
H(x,p):=F(x,p)+G(x). \label{H}%
\end{equation}
Using Lemma \ref{lip_H_par}, we know that $H$ is Lipschitz-like with respect
to $x$ uniformly in $p$ around $((\overline{x},\overline{p}),{0})$ and the
relation (\ref{lip_sum_par}) holds for $\overline{z}:=-{\overline{y}.}$

Using now Proposition \ref{link_around}, (iii) is equivalent to the fact that
$F$ is open at linear rate with respect to $p$ uniformly in $x$ around
$((\overline{x},\overline{p}),\overline{y})$ and $\widehat{\operatorname*{lop}%
}_{p}F((\overline{x},\overline{p}),\overline{y})=(\widehat{\operatorname*{reg}%
}_{p}F((\overline{x},\overline{p}),\overline{y}))^{-1}.$ Consequently, there
exist $\varepsilon,L>0$ such that, for every $x\in B({\overline{x}%
,}\varepsilon)$, every $(p,y)\in\operatorname*{Gr}F_{x}\cap\lbrack
B(\overline{p},\varepsilon)\times B(\overline{y},\varepsilon)]$ and every
$\rho\in]0,\varepsilon\lbrack,$ one has%
\[
B(y,L\rho)\subset F_{x}(B(p,\rho)).
\]

Using now (i), there exists $\delta\in]0,\varepsilon\lbrack$ such that, for
every $(x,p)\in B(\overline{x},\delta)\times B(\overline{p},\delta)$ and every
$w\in(F_{p}+G)(x)\cap B({0},\delta),$ there exist $y\in F_{p}(x)\cap
B({\overline{y},}\varepsilon)$ and $z\in G(x)\cap B(-\overline{y}%
,\varepsilon)$ such that $w=y+z.$

Take now arbitary $x\in B({\overline{x},}\delta)$, $(p,w)\in\operatorname*{Gr}%
H_{x}\cap\lbrack B(\overline{p},\delta)\times B(0,\delta)]$ and $\rho
\in]0,\varepsilon\lbrack.$ Then $w\in H(x,p)\cap B({0},\delta),$ so $w$ can be
written as $y+z,$ with $(p,y)\in\operatorname*{Gr}F_{x}\cap\lbrack
B(\overline{p},\varepsilon)\times B(\overline{y},\varepsilon)]$ and $z\in
G(x)\cap B(-\overline{y},\varepsilon),$ whence%
\[
B(w,L\rho)=z+B(y,L\rho)\subset G(x)+F(x,B(p,\rho))=H_{x}(B(p,\rho)).
\]

Adding the fact that $L$ can be chosen arbitrarly close to $\widehat
{\operatorname*{lop}}_{p}F((\overline{x},\overline{p}),\overline{y}),$ we
obtain that $H$ is open at linear rate with respect to $p$ uniformly in $x$
around $((\overline{x},\overline{p}),0)$ and $\widehat{\operatorname*{lop}%
}_{p}H((\overline{x},\overline{p}),0)=(\widehat{\operatorname*{reg}}%
_{p}F((\overline{x},\overline{p}),\overline{y}))^{-1}.$

Now the result easily follows from the second item of Theorem \ref{impl}%
.$\hfill\square$

\bigskip

Theorem \ref{mreg_sol} is a generalization of \cite[Theorem 3.3 (i)]%
{ArtMord2009}, concerning the direct implication.\textbf{ }

Notice that in \cite[Theorems 3.2, 3.3]{Uderzo2009} some related results are
obtained as well. However, Uderzo's approach, which works on metric spaces,
concerns global openness results and involves a parametric function instead of
our multifunction $F$, and $G$ is also taken in parametric form.

On the other hand, the next result uses the ideas of the proof of Theorem
\ref{impl} and, essentially, the estimations from Theorem \ref{main_const}.
For several technical reasons it is not possible to give a direct and easy
proof based on the main results, but, nevertheless, the proof uses the very
same ideas and arguments arranged in the specific context of this result.

Let us emphasize that, once again, the Lipschitz modulus of $S$ has a form
which can be easily developed from the previous facts. Namely, one can expect
that $\operatorname*{lip}S(\overline{p},\overline{x})\leq c^{-1}%
\widehat{\operatorname*{lip}}_{p}H((\overline{x},\overline{p}),0)=c^{-1}%
\widehat{\operatorname*{lip}}_{p}F((\overline{x},\overline{p}),\overline{y}),$
where $c>0$ must be the rate of linear openness of $H$ with respect to $x.$
One can observe that for every $p$ in a neighborhood of $\overline{p},$ and
for appropriate $x,y,$ $\widehat{\operatorname*{lip}}_{x}F((\overline
{x},\overline{p}),\overline{y})$ seems to be very close to
$\operatorname*{lip}(F_{p})(x,y)=$ $\operatorname*{lip}(-F_{p})(x,-y)=$
$\left(  \operatorname*{lop}(-F_{p})^{-1}(-y,x)\right)  ^{-1}.$ Also,
$\operatorname*{reg}G(\overline{x},-{\overline{y}})=$ $\left(
\operatorname*{lop}G(\overline{x},-{\overline{y}})\right)  ^{-1}.$ Because we
know additionally that $\widehat{\operatorname*{lip}}_{x}F((\overline
{x},\overline{p}),\overline{y})\cdot\operatorname*{reg}G(\overline
{x},-{\overline{y}})<1,$ then one can expect that $G-(-F_{p})$ might be open
at the linear rate $c=\operatorname*{lop}G(\overline{x},-{\overline{y}%
})-\widehat{\operatorname*{lip}}_{x}F((\overline{x},\overline{p}),\overline
{y})$ $=\left(  \operatorname*{reg}G(\overline{x},-{\overline{y}})\right)
^{-1}-\widehat{\operatorname*{lip}}_{x}F((\overline{x},\overline{p}%
),\overline{y})$ at points close to $((\overline{x},\overline{p}),0).$ Whence,
once again, one can have an intuitive approach in getting (\ref{lS}) below.

\begin{thm}
\label{Lip_sol}Let $X,Y,P$ be Banach spaces, $F:X\times P\rightrightarrows Y,$
$G:X\rightrightarrows Y$ be two set-valued maps and $(\overline{x}%
,\overline{p},\overline{y})\in X\times P\times Y$ such that $\overline{y}\in
F(\overline{x},\overline{p})$ and $-\overline{y}\in G(\overline{x})$. Suppose
that the following assumptions are satisfied:

(i) $(F,G)$ is locally sum-stable with respect to $x$ uniformly in $p$ around
$(\overline{x},\overline{p},{\overline{y},}-{\overline{y}});$

(ii) for every $p$ in a neighborhood of $\overline{p},$ $\operatorname*{Gr}%
F_{p}$ is closed;

(iii) $\operatorname{Gr}G$ is closed;

(iv) $F$ is Lipschitz-like around $((\overline{x},\overline{p}),\overline
{y});$

(v) $G$ is metrically regular around $(\overline{x},-\overline{y});$

(vi) $\widehat{\operatorname*{lip}}_{x}F((\overline{x},\overline{p}%
),\overline{y})\cdot\operatorname*{reg}G(\overline{x},-{\overline{y}})<1.$

Then $S$ is Lipschitz-like around $(\overline{p},\overline{x})$. Moreover, the
next relation is satisfied%
\begin{equation}
\operatorname*{lip}S(\overline{p},\overline{x})\leq\frac{\operatorname*{reg}%
G(\overline{x},-{\overline{y}})\cdot\widehat{\operatorname*{lip}}%
_{p}F((\overline{x},\overline{p}),\overline{y})}{1-\widehat
{\operatorname*{lip}}_{x}F((\overline{x},\overline{p}),\overline{y}%
)\cdot\operatorname*{reg}G(\overline{x},-{\overline{y}})}. \label{lS}%
\end{equation}

\end{thm}

\noindent\textbf{Proof. }Take $m>\operatorname*{reg}G(\overline{x}%
,-{\overline{y}})$ and $l>\widehat{\operatorname*{lip}}_{x}F((\overline
{x},\overline{p}),\overline{y})$ such that $m\cdot l<1.$

We intend to prove that there exist $\tau,t,\gamma>0$ such that, for every
$(x,p)\in B(\overline{x},\tau)\times B(\overline{p},t),$%
\begin{equation}
d(x,S(p))\leq m(1-lm)^{-1}d(0,[F(x,p)+G(x)]\cap B(0,\gamma)). \label{mrg}%
\end{equation}

Using the assumption (iv), one can find $r_{1},t_{1}>0$\ such that, for every
$p\in B(\overline{p},t_{1}),$ and every $x,x^{\prime}\in B(\overline{x}%
,r_{1}),$%
\begin{equation}
F_{p}(x)\cap B(\overline{y},s_{1})\subset F_{p}(x^{\prime})+l\left\Vert
x-x^{\prime}\right\Vert \mathbb{D}_{Y}. \label{lipFp}%
\end{equation}

But this shows, as one can see next, that for every $p\in B(\overline{p}%
,t_{1}),$ $F_{p}$ is $l-$pseudocalm at every $(x,y)\in\operatorname*{Gr}%
F_{p}\cap\lbrack B(\overline{x},2^{-1}r_{1})\times B(\overline{y},s_{1})].$
Indeed, take $x\in B(\overline{x},2^{-1}r_{1})$ and $y\in F_{p}(x)\cap
B(\overline{y},s_{1}).$ Then for every $x^{\prime}\in B(x,2^{-1}r_{1}),$ we
have from (\ref{lipFp}) that there exists $y^{\prime}\in F_{p}(x^{\prime})$
such that%
\[
d(y,F_{p}(x^{\prime}))\leq\left\Vert y-y^{\prime}\right\Vert \leq l\left\Vert
x-x^{\prime}\right\Vert ,
\]

\noindent which proves the desired assertion. Hence, we conclude in view of
Proposition \ref{link_at} that for every $p\in B(\overline{p},t_{1}),$
$F_{p}^{-1}$ is $l^{-1}-$open at every $(y,x)\in\operatorname*{Gr}F_{p}%
^{-1}\cap\lbrack B(\overline{y},s_{1})\times B(\overline{x},2^{-1}r_{1})].$

Also, from (v), we know that there exist $r_{2},s_{2}>0$ such that, for every
$(u,v)\in\operatorname{Gr}G\cap\lbrack B(\overline{x},r_{2})\times
B(-\overline{y},s_{2})],$ $G$ is metrically hemiregular at $(u,v)$ with
constant $m,$ whence is open at linear rate $m^{-1}$ at $(u,v).$

Use now the property from (i) for $\min\{2^{-1}s_{1},2^{-1}s_{2}\}$ instead of
$\varepsilon$ and find $\delta$ such that the assertion from Definition
\ref{sum-sta_param} is true.

Take now $\rho\in]0,\min\{(m^{-1}-l)^{-1}\delta,4^{-1}r_{1},2^{-1}r_{2}%
,2^{-1}l^{-1}s_{1},2^{-1}ms_{2}\}[$ and define $\gamma:=(m^{-1}-l)\rho.$
Suppose that $B(\overline{p},t_{2})$ is the neighborhood which appears in
(ii), denote by $t:=\min\{t_{1},t_{2}\},$ $\tau:=\min\{4^{-1}r_{1},2^{-1}%
r_{2}\},$ and take $(x,p)\in B(\overline{x},\tau)\times B(\overline{p},t).$

If $[F(x,p)+G(x)]\cap B(0,\gamma)=\emptyset$ or $0\in\lbrack F(x,p)+G(x)]\cap
B(0,\gamma),$ then (\ref{mrg}) trivially holds. Suppose that $0\not \in
\lbrack F(x,p)+G(x)]\cap B(0,\gamma).$ Then, for every $\varepsilon>0,$ one
can find $w_{\varepsilon}\in\lbrack F(x,p)+G(x)]\cap B(0,\gamma)$ such that%
\begin{equation}
\left\Vert w_{\varepsilon}\right\Vert <d(0,[F(x,p)+G(x)]\cap B(0,\gamma
))+\varepsilon. \label{ineg}%
\end{equation}

Obviously, $d(0,[F(x,p)+G(x)]\cap B(0,\gamma))<(m^{-1}-l)\rho,$ so for
$\varepsilon>0$ sufficiently small, $d(0,[F(x,p)+G(x)]\cap B(0,\gamma
))+\varepsilon<(m^{-1}-l)\rho.$ Consequently, we get from (\ref{ineg}) that%
\begin{equation}
0\in B(w_{\varepsilon},d(0,[F(x,p)+G(x)]\cap B(0,\gamma))+\varepsilon)\subset
B(w_{\varepsilon},(m^{-1}-l)\rho)=B(w_{\varepsilon},\delta), \label{0in}%
\end{equation}

\noindent Applying (i), one can find $y_{\varepsilon}\in F(p,x)\cap
B({\overline{y},2}^{-1}s_{1})$ and $z_{\varepsilon}\in G(x)\cap B(-{\overline
{y},2}^{-1}s_{2})$ such that $w_{\varepsilon}=y_{\varepsilon}+z_{\varepsilon
}.$ Whence, $B(y_{\varepsilon},2^{-1}s_{1})\subset B(\overline{y},s_{1})$ and
$B(z_{\varepsilon},2^{-1}s_{2})\subset B(-\overline{y},s_{2}).$

Denote now $r_{1}^{\prime}:=4^{-1}r_{1},$ $s_{1}^{\prime}:=2^{-1}s_{1},$
$r_{2}^{\prime}:=2^{-1}r_{2},$ $s_{2}^{\prime}:=2^{-1}s_{2}.$ Summarizing,
$F_{p}^{-1}$ is $l^{-1}-$open at every $(y^{\prime},x^{\prime})\in
\operatorname*{Gr}F_{p}^{-1}\cap\lbrack B(y_{\varepsilon},s_{1}^{\prime
})\times B(x,r_{1}^{\prime})]$, $G$ is $m^{-1}-$ open at every $(u^{\prime
},v^{\prime})\in\operatorname*{Gr}G\cap\lbrack B(x,r_{2}^{\prime})\times
B(z_{\varepsilon},s_{2}^{\prime})]$ and $l^{-1}m^{-1}>1.$ We can apply then
Theorem \ref{main_const} for $-G,$ $F_{p}^{-1},$ $(x,-z_{\varepsilon}%
)\in\operatorname*{Gr}(-G),$ $(y_{\varepsilon},x)\in\operatorname*{Gr}%
F_{p}^{-1}$ and $\rho_{0}:=(m^{-1}-l)^{-1}(d(0,[F(x,p)+G(x)]\cap
B(0,\gamma))+\varepsilon)<\rho<\min\{r_{1}^{\prime},r_{2}^{\prime},l^{-1}%
s_{1}^{\prime},ms_{2}^{\prime}\}$ to obtain that%
\[
B(y_{\varepsilon}+z_{\varepsilon},d(0,[F(x,p)+G(x)]\cap B(0,\delta
))+\varepsilon)\subset(F_{p}+G)(B(x,\rho_{0})).
\]

Using (\ref{0in}), we obtain that $0\in(F_{p}+G)(B(x,\rho_{0})),$ so there
exists $\widetilde{x}\in B(x,\rho_{0})$ such that $0\in F(\widetilde
{x},p)+G(\widetilde{x})$ or, equivalently, $\widetilde{x}\in S(p).$ Hence%
\[
d(x,S(p))\leq\left\Vert x-\widetilde{x}\right\Vert <\rho_{0}=(m^{-1}%
-l)^{-1}d(0,[F(x,p)+G(x)]\cap B(0,\delta))+\varepsilon.
\]

Making $\varepsilon\rightarrow0,$ we obtain (\ref{mrg}).

For the final step of the proof, observe that because $F$ is Lipschitz-like
with respect to $p$ uniformly in $x$ around $((\overline{x},\overline
{p}),{\overline{y}}),$ there exist $\alpha,k>0$ such that for every $x\in
B({\overline{x},}\alpha),$ every $p_{1},p_{2}\in B(\overline{p},\alpha),$ one
has%
\begin{equation}
F_{x}(p_{1})\cap B(\overline{y},\alpha)\subset F_{x}(p_{2})+k\left\Vert
p_{1}-p_{2}\right\Vert \mathbb{D}_{Y}. \label{lipFx}%
\end{equation}

Use now (i) for $\alpha$ instead of $\varepsilon$ and find $\delta^{\prime}%
\in]0,\alpha\lbrack$ such that the assertion from Definition
\ref{sum-sta_param} is true. Take now arbitrary $x\in B({\overline{x},}%
\delta^{\prime}),$ $p_{1},p_{2}\in B(\overline{p},\delta^{\prime})$ and $w\in
H_{x}(p_{1})\cap B(0,\delta^{\prime}),$ where $H$ is defined by (\ref{H}).
Thene there exist $y\in F_{x}(p_{1})\cap B(\overline{y},\alpha)$ and $z\in
G(x)\cap B(-\overline{y},\alpha)$ such that $w=y+z.$ Using now (\ref{lipFx}),
we have that%
\begin{align*}
w  &  =y+z\in F(x,p_{2})+k\left\Vert p_{1}-p_{2}\right\Vert \mathbb{D}%
_{Y}+G(x)\\
&  =H_{x}(p_{2})+k\left\Vert p_{1}-p_{2}\right\Vert \mathbb{D}_{Y}.
\end{align*}

\noindent Consequently, $H$ is Lipschitz-like with respect to $p$ uniformly in
$x$ around $((\overline{x},\overline{p}),{0})$. Moreover, because $k$ can be
chosen arbitrarly close to $\widehat{\operatorname*{lip}}_{p}F((\overline
{x},\overline{p}),\overline{y}),$ we have that $\widehat{\operatorname*{lip}%
}_{p}F((\overline{x},\overline{p}),\overline{y})=\widehat{\operatorname*{lip}%
}_{p}H((\overline{x},\overline{p}),0).$ Then one can proceed as in the proof
of the final part of Theorem \ref{impl}\ (i) to conclude the proof.$\hfill
\square$

\bigskip

The result would follow using just the local closedness of the graphs of
$F_{p}$ and $G$ for $p$ close to $\overline{p}$ around appropriate points, but
we preferred the actual formulation in order to avoid a more complicated
assertion. Theorem \ref{Lip_sol} generalizes \cite[Theorem 5.1 (ii)]%
{ArtMord2009}.

\bigskip

Note that, using Proposition \ref{gama}, one could obtain on the basis of
Artacho and Mordukhovich tools in \cite{ArtMord2009} the converse links
between the metric regularity/Lipschitz-like property of solution map $S$ and
the corresponding properties of the field map $G.$

\bigskip

\noindent{\small {\textbf{Acknowledgement: }}We are indebted to an anonymous
referee for several constructing remarks which led to the improvement of our
work.}

\end{document}